\input amstex\input epsf
\documentstyle{amsppt}\nologo\footline={}

\def\B{\mathop{\Bbb B}}
\def\Isom{\mathop{\text{\rm Isom}}}
\def\dist{\mathop{\text{\rm dist}}}
\def\PU{\mathop{\text{\rm PU}}}
\def\Lin{\mathop{\text{\rm Lin}}}
\def\tr{\mathop{\text{\rm tr}}}
\def\S{\mathop{\Bbb S}}
\def\GL{\mathop{\text{\rm GL}}}
\def\Re{\mathop{\text{\rm Re}}}
\def\Im{\mathop{\text{\rm Im}}}
\def\T{\mathop{\text{\rm T}}}
\def\Stab{\mathop{\text{\rm Stab}}}
\def\U{\mathop{\text{\rm U}}}
\def\SU{\mathop{\text{\rm SU}}}
\def\o{\mathop{\text{\rm o}}}
\def\O{\mathop{\text{\rm O}}}
\def\Gr{\mathop{\text{\rm Gr}}}
\def\ssl{\mathop{\text{\rm sl}}}
\def\Sym{\mathop{\text{\rm Sym}}}

\hsize450pt\topmatter\title Complex hyperbolic equidistant loci\endtitle\author Sasha Anan$'$in\endauthor\address Departamento de
Matem\'atica, ICMC, Universidade de S\~ao Paulo, Caixa Postal 668, \newline13560-970--S\~ao Carlos--SP,
Brasil\endaddress\subjclass51M10 (57S30, 17A99)\endsubjclass\abstract We describe and study the loci equidistant from finitely many
points in the so-called complex hyperbolic geometry, i.e., in the geometry of a holomorphic $2$-ball $\B$. In particular, we show that
the bisectors (= the loci equidistant from $2$ points) containing the (smooth real algebraic) curve equidistant from given $4$ generic
points form a real elliptic curve and that the foci of the mentioned bisectors constitute an isomorphic elliptic curve.

We are going to use the obtained facts in constructions of (compact) quotients of $\B$ by discrete groups.

With similar technique, we also classify up to isotopy generic $3$-dimensional algebras (i.e., bilinear operations) over an
algebraically closed field $\Bbb K$ of characteristic $\ne2,3$. Briefly speaking, an algebra is classified by the (plane projective)
curve $D$ of its zero divisors equipped with a nonprojective automorphism of $D$. This classification is almost equivalent to the
classification of the so-called geometric tensors given in [BoP] by A.~Bondal and A.~Polishchuk in their study of noncummutative
projective planes.\endabstract\endtopmatter\document

{\hfill\it To Victor Gerasimov on occasion of\/ $50$ years of}

{\hfill\it our mathematical discussions}

\bigskip

\bigskip

\centerline{\bf1.~Introduction}

\medskip

A group $G\subset\Isom X$ of isometries of a simply connected (model) space $X$ is known to be discrete iff a point $p\in X$ not fixed
by any nontrivial element of $G$ lies in the interior of the {\it Dirichlet polyhedron\/}
$P:=\{x\in X\mid\dist(x,p)\le\dist(x,gp)\text{ for all }g\in G\}$ centred at $p$. The faces of $P$ are frequently loci equidistant
from finitely many points $p,g_1p,\dots,g_np$.

Dealing with the geometry of the holomorphic $2$-ball $\B$, it is possible to almost recover (up to $4$ real parameters corresponding
to the freedom in the choice of $p$) a cocompact discrete group $G\subset\Isom\B=\PU(2,1)$ from the combinatorics of the face pairing
of its Dirichlet polyhedron [Ana]. Therefore, wishing to construct such groups, it is important to study complex hyperbolic
equidistant loci.

The well-known G.~Giraud rigidity theorem [Gol, Theorem 8.3.3, p.~264] establishes that there are at most $3$ bisectors containing the
nonempty locus equidistant from $3$ noncollinear points. (By definition, a~bisector is the locus equidistant from $2$ points.) As
noted in [Gol, pp.~viii, x, xiv], this theorem constitutes the main constraint on the combinatorics of Dirichlet polyhedra.
Indeed, it is possible to reduce the study of discrete cocompact groups to the case where, in the tessellation of $\B$ by the
polyhedra congruent to $P$, every codimension $2$ face is contained in exactly $3$ codimension $1$ faces [Ana]. In particular, every
defining relation between the face pairing isometries has length $3$. Moreover, almost all conditions of Poincar\'e's polyhedron
theorem [AGS] follow from the Giraud rigidity. In our exposition (see Lemma 4.7 and Corollary 4.8), the Giraud theorem sounds almost
as `there are at most $3$ roots of a polynomial of degree $3$'. In spite of this, our results can be seen as a development of the
mentioned Giraud rigidity.

\medskip

{\bf1.1.~Bisectors and their ingredients.} We follow the notation in [AGr] and [AGG] (for a background in complex hyperbolic geometry,
see also [Gol]).

Throughout the paper, $V$ is a $3$-dimensional $\Bbb C$-linear space equipped with a hermitian form $\langle-,-\rangle$ of signature
$++-$ so that all negative points in the complex projective plane constitute the holomorphic $2$-ball
$\B:=\big\{p\in\Bbb P_\Bbb CV\mid\langle p,p\rangle<0\big\}$. The distance on $\B$ is given by
$\cosh^2\dist(p_1,p_2)=\frac{\langle p_1,p_2\rangle\langle p_2,p_1\rangle}{\langle p_1,p_1\rangle\langle p_2,p_2\rangle}$.

Hence, the locus equi(dis)tant from $2$ distinct points $p_1,p_2\in\B$ is given by the equation $[x,x]=0$, where the hermitian form
$[x,y]:=\frac{\langle x,p_1\rangle\langle p_1,y\rangle}{\langle p_1,p_1\rangle}-\frac{\langle x,p_2\rangle\langle
p_2,y\rangle}{\langle p_2,p_2\rangle}$
has rank $2$. Since $[-,-]=\langle h-,-\rangle$ for a suitable (unique) $\Bbb C$-linear map $h\in\Lin_\Bbb C(V,V)$ of rank $2$ such
that $h^*=h$ and $\tr h=0$, where $h^*$ denotes the map adjoint to $h$ in the sense of $\langle-,-\rangle$, we arrive at the
following definition.

\medskip

{\bf1.1.1.~Definition.} Let $h\in\Lin_\Bbb C(V,V)$ be a $\Bbb C$-linear map of rank $2$ such that $h^*=h$ and $\tr h=0$. Then
$B_h:=\big\{p\in\Bbb P_\Bbb CV\mid\langle hp,p\rangle=0\big\}$ is the {\it bisector\/} given by $h$ if $B_h\ne\{f\}$, where
$\{f\}:=\Bbb P_\Bbb C\ker h$ is the {\it focus\/} of $B_h$. The projective line $C:=\Bbb P_\Bbb ChV=\Bbb P_\Bbb Cf^\perp$ is the {\it
complex spine\/} of $B_h$. Obviously, $B_h$ is a projective cone with apex $f$. The projective line $S_p\subset B_h$ spanned by $f$
and $p\in B_h\setminus\{f\}$ is the {\it slice\/} $S_p$ of $B_h$ generated by $p$. The point $hp$ is the point polar (= orthogonal) to
the slice $S_p$. The points polar to the slices form a geodesic $R\subset C$ called the {\it real spine\/} of $B_h$. The bisector
$B_h$ is {\it hyperbolic,spherical,} or {\it parabolic\/} if its focus $f$ is respectively positive (i.e., $\langle f,f\rangle>0$),
negative (i.e., $f\in\B$), or isotropic (i.e., $f\in\S:=\big\{p\in\Bbb P_\Bbb CV\mid\langle p,p\rangle=0\big\}$).

\medskip

{\bf1.1.2.~Remark.} {\sl Bisectors\/ $B_{h_1},B_{h_2}$ coincide iff\/ $h_1,h_2$ are\/ $\Bbb R^*$-proportional.}

\medskip

Our definition of a bisector differs from the commonly accepted one; the latter deals with the bisectors which are hyperbolic in our
sense. It follows a motivation of our definition. Although a spherical bisector is singular at its focus, parts of such bisectors can
still constitute smooth codimension $1$ faces of fundamental polyhedra. Thus, we can work with Dirichlet polyhedra centred at positive
points. The other reason is that spherical and parabolic bisectors appear naturally when we deform the usual hyperbolic ones;
for~example, the elliptic family of bisectors in Theorem 1.2.11 can easily contain all $3$ types.

\medskip

{\bf1.2.~Intersection and families of bisectors.} In order to motivate the concept of a family of bisectors, we begin with the
following easy but useful criterion.

\medskip

{\bf1.2.1.~Remark.} {\sl The bisectors\/ $B_{h_1},\dots,B_{h_n}$ are transversal at\/ $p\in\bigcap_iB_{h_i}\setminus\S$ iff\/
$\dim_\Bbb RW=n$ and\/ $hp\ne0$ for any $0\ne h\in W$, where\/ $W$ stands for the\/ $\Bbb R$-span of\/ $h_1,\dots,h_n$.}

\medskip

Note that the intersection in Remark 1.2.1 can be described as
$\bigcap_iB_{h_i}=\big\{p\in\Bbb P_\Bbb CV\mid\langle Wp,p\rangle=0\big\}$.

\medskip

{\bf1.2.2.~Definition.} An $\Bbb R$-linear subspace $W\subset\Lin_\Bbb C(V,V)$ is called a {\it family\/} of bisectors if $\tr W=0$,
$B_W\cap\B\ne\varnothing$, $W$ is the $\Bbb R$-span of $D_W:=\{h\in W\mid\det h=0\}$, and $h^*=h$ for all $h\in W$, where
$B_W:=\big\{p\in\Bbb P_\Bbb CV\mid\langle Wp,p\rangle=0\big\}$ is the {\it base\/} of the family. It is easy to see (Remark 3.4) that
any $0\ne h\in D_W$ has rank $2$, i.e., $B_h$ is a bisector. The image $E_W$ of $D_W$ in the real projective space $\Bbb P_\Bbb RW$ is
given by a single (possibly trivial) equation $\det=0$. If $D_W=W$, the family is {\it linear.} If $D_W\ne W$ and $\dim_\Bbb RW=3$,
the~family is {\it elliptic\/}; in this case $E_W$ is a real cubic; denote by $\hat E_W\subset\Bbb P_\Bbb C\Bbb CW$ the corresponding
complex plane cubic.

The algebraic map $f:E_W\to\Bbb P_\Bbb CV$ sending a bisector to its focus is called {\it focal.}

\medskip

The following proposition claims that, excluding the relatively trivial case of a confocal linear family, there are two types of
linear families: 1.~$E_W$ is a real projective plane of bisectors isomorphic by means of the focal map to an $\Bbb R$-plane;
2.~$E_W$ is a real projective space of bisectors sharing a common slice.

\medskip

{\bf1.2.3.~Proposition.} {\sl Let\/ $W$ be a\/ {\rm nonconfocal} {\rm(}i.e., $Wf_0\ne 0$ for any\/ $0\ne f_0\in V${\rm)} linear
family of bisectors. Then there exists a unique up to\/ $\Bbb C^*$-proportionality\/ $\Bbb R$-linear embedding
$f:W\hookrightarrow V$ such that\/ $h(fh)=0$ for all\/ $h\in W$. Denote\/ $U:=\Bbb CfW$. If\/ $U=V$, then\/ $\dim_\Bbb RW=3$ and\/
$fW\subset V$ is a totally real subspace. Otherwise, $\dim_\Bbb CU=2$ and the bisectors of the family share the common slice\/
$\Bbb P_\Bbb CU$, i.e., $W\subset W_U:=\big\{h\in\Lin_\Bbb C(V,V)\mid\tr h=0,\ h^*=h,\ hU\subset U^\perp\big\}$.}

\medskip

{\bf1.2.4.~Proposition.} {\sl With the exception of a confocal\/ $3$- or\/ $4$-dimensional linear family of bisectors with negative
common focus, the base\/ $B_W$ completely determines a linear family\/ $W$ because, in this case,
$W=\big\{h\in\Lin_\Bbb C(V,V)\mid\tr h=0,\ h^*=h,\ \langle hb,b\rangle=0\text{ \rm for all }b\in B_W\big\}$.}

\medskip

{\bf1.2.5.~Definition.} An elliptic family $W$ of bisectors is said to be {\it generic\/} if $E_W$ contains no confocal line and is
not a real projective line plus a point.

\medskip

{\bf1.2.6.~Lemma.} {\sl Let\/ $W$ be a generic elliptic family of bisectors. Then the focal map\/ $f$ can be extended to a\/ {\rm
focal} isomorphism\/ $\hat f:\hat E_W\to E\subset\Bbb P_\Bbb CV$ between complex plane cubics.}

\medskip

{\bf1.2.7.~Definition.} A generic elliptic family $W$ of bisectors is {\it real\/} if there are $3$ distinct points in $E_W$
whose foci lie on a same complex projective line not included in $E$.

\medskip

{\bf1.2.8.~The elliptic family equitant from $4$ points.} Let $p_i\in\Bbb P_\Bbb CV$, $1\le i\le4$, be points of a same signature
$\sigma$ such that no $3$ of them are on a same complex projective line. Then $w_i:=\langle-,p_i\rangle p_i$, $1\le i\le4$, are
$\Bbb C$-linearly independent due to an essentially unique $\Bbb C$-linear dependence between the~$p_i$'s. We can pick
representatives such that $\langle p_i,p_i\rangle=\sigma$. Suppose that there exists $p\in\B$ such that
$\big|\langle p,p_i\rangle\big|=1$ for all $1\le i\le 4$. Then
$W:=\big\{\sum_{i=1}^4r_iw_i\mid r_i\in\Bbb R,\ \sum_{i=1}^4r_i=0\big\}$ is an elliptic family of bisectors because
$\sum_{i=1}^4r_iw_i=\sum_{i=1}^3\langle-,p_i\rangle r_ip_i$ has rank $3$ if $r_4=\sum_{i=1}^3r_i=0$ and $r_i\ne0$ for all $i=1,2,3$.

The constructed elliptic family $W$ is called {\it equitant\/} from the $p_i$'s. When $\sigma=0$, the family $W$ depends on the
choice of representatives of the $p_i$'s.

By Lemma 5.4, if the curve $E_W$ of an equitant family $W$ is not irreducible, then it is a smooth real conic plus a real projective
line intersecting the conic in $2$ points or it consists of $3$ real projective lines sharing no common point.

\medskip

{\bf1.2.9.~Proposition.} {\sl An elliptic family of bisectors is equitant iff it is real.}

\medskip

{\bf1.2.10.~Proposition.} {\sl Let $W$ be an equitant elliptic family of bisectors with irreducible $E_W$. Then the base $B_W$ is an
irreducible algebraic curve.}

\medskip

The main theorem of the paper says that the bisectors containing an infinite subset of the locus equi(dis)tant from given $4$ generic
points of a same signature form a real elliptic curve. Under a certain angle of view, this fact is analogous to the Giraud rigidity
theorem. Note that the elliptic curve can easily have bisectors of all $3$ types.

\medskip

{\bf1.2.11.~Theorem.} {\sl Let\/ $W$ be an equitant elliptic family of bisectors with irreducible\/ $E_W$ and let\/ $B\subset B_W$ be
infinite. Then\/ $W=\big\{h\in\Lin_\Bbb C(V,V)\mid h^*=h,\ \langle hb,b\rangle=0\text{ \rm for all }b\in B\big\}$.}

\medskip

{\bf1.3.~Classification of $3$-dimensional algebras up to isotopy.} In fact, this theme is not related to the previous one, and the
reader interested only in equidistant loci may simply ignore it. The one who is interested in proofs can find them in the appendix to
this paper.

Let $\Bbb K$ be a field and let $V_i$, $i=1,2,3$, be $3$-dimensional $\Bbb K$-linear spaces. We consider the task of classifying
elements in $V_1\otimes_\Bbb KV_2\otimes_\Bbb KV_3$ modulo the action of the group $\GL V_1\times\GL V_2\times\GL V_3$.

The addressed classification is equivalent to the classification of $3$-dimensional $\Bbb K$-algebras, i.e.,\break
$\Bbb K$-bilinear maps $b:A\times A\to A$, $\dim_\Bbb KA=3$, modulo the action of the group $\GL A\times\GL A\times\GL A$ given by the
rule $b^{(g_1,g_2,g_3)}(a_1,a_2):=g_3b(g_1^{-1}a_1,g_2^{-1}a_2)$ for all $a_1,a_2\in A$ and $g_1,g_2,g_3\in\GL A$. Algebras in a same
orbit are said to be {\it isotopic.} The detailed classification of $3$-dimensional algebras modulo isotopy is too bulky and boring to
be described here, especially if $\Bbb K$ is not algebraically closed. This is why we give the answer only for generic $3$-dimensional
algebras over an algebraically closed $\Bbb K$. We denote the bilinear operation $b$ by $\cdot$.

When $d_1\cdot d_2=0$ and $0\ne d_1,d_2\in A$, we call $d_1$ and $d_2$ nontrivial left and right zero divisors. Usually, we regard
them as points in the projective plane $\Bbb P_\Bbb KA$. Denote by $D_1,D_2\subset\Bbb P_\Bbb KA$ the schemes of left/right zero
divisors. The scheme $D_1$ is given in $\Bbb P_\Bbb KA$ by the equation $p(x_0,x_1,x_2)=0$, where $p(x_0,x_1,x_2):=\det\Phi$ is a
homogeneous polynomial of degree $3$ (or $p$ equals $0$), $\Phi:=a_0x_0+a_1x_1+a_2x_2$, the elements of $A$ are considered as left
multiplications, and $a_0,a_1,a_2$ is a $\Bbb K$-linear basis in $A$.

\medskip

{\bf1.3.1.~Definition.} We call a $3$-dimensional $\Bbb K$-algebra $A$ {\it generic\/} if the elements in $D_1$ and $D_2$ considered
respectively as left and right multiplications have rank $2$. For a generic $A$, we obtain an isomorphism $\varphi:D_1\to D_2$ given
by the relation $d_1\cdot d_2=0$. An algebra isotopic to $A$ provides the isomorphism $i_2\varphi i_1:D_1\to D_2$, where $i_j$ is a
projective automorphism of $D_j$.

\medskip

If the scheme $D_1$ of left zero divisors is a smooth cubic, the $3$-dimensional algebra is necessarily generic because $D_1$ is not
rational.\footnote{For any other $D_1$, i.e., for any non-smooth cubic $D_1$ (or for $D_1=\Bbb P_\Bbb KA$), there exists a nongeneric
algebra.}
(The following obvious remark is silently used here: If $a_i\cdot L_i=0$, $i=1,2$, then $L\cdot p=0$, where $L_i\subset\Bbb P_\Bbb KA$
are lines, $L\subset\Bbb P_\Bbb KA$ is the line spanned by distinct $a_1,a_2\in\Bbb P_\Bbb KA$, and $p\in L_1\cap L_2$.)

\medskip

{\bf1.3.2.~Theorem.} {\sl Let\/ $\Bbb K$ be an algebraically closed field of characteristic\/ $\ne2,3$.

Up to isotopy, the\/ $3$-dimensional generic\/ $\Bbb K$-algebras\/ $A$ whose zero divisors variety\/ $D$ is isomorphic to a cubic
without multiple components are classified by a nonprojective automorphism of\/ $D$ {\rm(}such an automorphism always exists\/{\rm)}
considered modulo projective automorphisms.

When\/ $D$ is a double line plus a line or a triple line or\/ $D=\Bbb P_\Bbb KA$, there exists a unique, up to isotopy,
$3$-dimensional generic\/ $\Bbb K$-algebra\/ $A$ with the indicated\/ $D$.

Briefly speaking, a generic\/ $3$-dimensional\/ $\Bbb K$-algebra modulo isotopy is a plane cubic $D$ of its zero divisors equipped
with a nonprojective automorphism of $D$.}

\medskip

When studying noncommutative projective planes, A.~Bondal and A.~Polishchuk [BoP] classified the so-called geometric tensors. This
classification is almost equivalent to that of generic algebras (see [BoP, Table, p.~36] for details). The algebras whose variety of
zero divisors is a smooth conic plus a line such that the isomorphism between $D_1,D_2$ interchanges the conic and the line constitute
the difference between generic algebras and geometric tensors.

The complete classification of $4$-dimensional $\Bbb K$-algebras up to isotopy seems to be a difficult task. One can readily observe
that the smooth hypersurfaces in $\Bbb P_\Bbb K^3$ of degree $4$ (they are K3-surfaces) will occupy the place of smooth cubics. It is
curious that the scheme $D$ of zero divisors of the algebra of quaternions is the most `degenerate' K3-surface given by the equation
$(x_0^2+x_1^2+x_2^2+x_3^2)^2=0$.

\bigskip

\centerline{\bf2.~Basic properties of bisectors}

\medskip

Most of the following material is well known and can be found in [Gol] (sometimes, using different terminology and, usually, in a
different exposition).

\medskip

{\bf2.1.~Hyperbolic, spherical, and parabolic bisectors.} If $h$ possesses a nonnull eigenvalue $\lambda$, then, in view of $\tr h=0$,
we have the orthogonal $h$-invariant decomposition $V=\ker h\oplus hV$ and the eigenvectors $v_1,v_2$ of $h$ corresponding to
$\lambda,-\lambda$, called {\it vertices\/} of $B_h$, span the complex spine $C:=\Bbb P_\Bbb ChV$. The bisector $B_h$ is
hyperbolic/spherical iff $C$ has signature $+-$/$++$.

In both cases, $B_h$ is the projective cone with apex $f$ and base $C\cap B_h\ne\varnothing$. Moreover, $C\cap B_h$ is a
hyperbolic/spherical geodesic and coincides with the real spine $R$ of $B_h$.

Indeed, if $v_1,v_2$ are isotropic, then the eigenvalues $\lambda,-\lambda$ of $h$ are purely imaginary because
$0\ne\lambda\langle v_1,v_2\rangle=\langle hv_1,v_2\rangle=\langle v_1,hv_2\rangle=-\overline\lambda\langle v_1,v_2\rangle$. It is
easy to see that $C\cap B_h$ is the hyperbolic geodesic with vertices $v_1,v_2$.

If $v_1$ is nonisotropic, then $\lambda\in\Bbb R$ and $v_1,v_2$ are orthogonal. We choose orthonormal representatives of $v_1,v_2$.
When $v_1,v_2$ have a same signature, one can see that $C\cap B_h=\{\overline uv_1+uv_2\mid|u|=1\}$ is a spherical geodesic. When
$v_1,v_2$ have different signatures, the set $\big\{c\in C\mid\langle hc,c\rangle=0\big\}$ is empty; so, this case is impossible.
While $p$ runs over the geodesic $C\cap B_h$, the point polar to the slice $S_p$ runs over the geodesic $C\cap B_h$. This implies
$R=C\cap B_h$.

If the eigenvalues of $h$ are all null, then $h^3=0$. Since $\dim_\Bbb C\ker h=1$, the focus $\{f\}:=\Bbb P_\Bbb C\ker h$ of $B_h$
lies in the complex spine, $f\in C:=\Bbb P_\Bbb ChV$, implying that $C=\Bbb P_\Bbb Cf^\perp$ has signature $+0$ and that $f$ is
isotropic; hence, $B_h$ is parabolic. It follows from $\langle h^2p,hp\rangle=\langle h^3p,p\rangle=0$ that $C\subset B_h$.

Let $d\notin C=\Bbb P_\Bbb ChV$. Then $h^2d\ne0$. So, $h^2d=f$ and $\Bbb R\ni\langle hd,hd\rangle=\langle h^2d,d\rangle\ne0$. Taking
$p:=d+rhd$, we have $\langle hp,p\rangle=\langle hd,d\rangle+2r\langle h^2d,d\rangle=0$ for a suitable $r\in\Bbb R$. Since
$p\notin C$, we obtain a slice $S_p$ of $B_h$ different from $C$. As $S_p$ contains an isotropic point different from $f$, we can
assume that $f\ne p\in\S\cap B_h$. Let $S':=\Bbb P_\Bbb Cp^\perp$. Every slice $S_q$ of $B_h$ intersects $S'$ in exactly $1$ point and
this point belongs to $B_h$. In~other words, $B_h$ is the projective cone with apex $f$ and base $\Gamma:=S'\cap B_h$. As
$\langle h^2p,p\rangle=\langle hp,hp\rangle>0$ (the point $hp$ is orthogonal to the isotropic point $p$, hence, is positive), we can
choose a representative $p$ such that the $\Bbb C$-linear basis $p,hp,h^2p$ of $V$ has the Gram matrix
$\left[\smallmatrix0&0&1\\0&1&0\\1&0&0\endsmallmatrix\right]$. Since $h=\left[\smallmatrix0&0&0\\1&0&0\\0&1&0\endsmallmatrix\right]$
in this basis, one can easily see that the base $\Gamma$ of the projective cone $B_h$ is the euclidean geodesic $\Bbb P_\Bbb CW'$,
where $W'$ is the $\Bbb R$-span of $p,ihp$. (There is no canonical choice for the geodesic $\Gamma$; see also the proofs of
Lemmas~4.4.3 and 4.4.4.) The real spine $R$ of $B_h$ is an euclidean geodesic as well: $R=\Bbb P_\Bbb CW$, where $W:=hW'$ is the
$\Bbb R$-span of $hp,ih^2p$.

\medskip

{\bf2.2.~Remark.} {\sl Let\/ $B_h$ be a parabolic bisector. Then there exists an isotropic point\/ $f\ne p\in\S\cap B_h$ different
from the focus $f$ of\/ $B_h$. For any such\/ $p\in B_h$, the points\/ $p,hp,h^2p$ form a\/ $\Bbb C$-linear basis in\/ $V$, this
basis has the Gram matrix\/ $\left[\smallmatrix0&0&1\\0&1&0\\1&0&0\endsmallmatrix\right]$ for a suitable choice of a representative\/
$p\in V$, and\/ $f=h^2p$. The bisector\/ $B_h$ is the projective cone with apex\/ $f$ and base\/ $\Gamma:=\Bbb P_\Bbb CW'=S'\cap B_h$,
where\/ $W'$ denotes the\/ $\Bbb R$-span of\/ $p,ihp$, $S':=\Bbb P_\Bbb Cp^\perp$, and\/ $\Gamma$ is an euclidean geodesic. The real
spine of\/ $B_h$ has the form\/ $R=\Bbb P_\Bbb CW$, where\/ $W\subset V$ stands for the\/ $\Bbb R$-span of\/ $hp,ih^2p$}
$_\blacksquare$

\medskip

Conversely, using the above considerations, it is easy to see that an arbitrary geodesic serves as a real spine of some bisector. Of
course, any bisector is uniquely determined by its real spine.

\medskip

{\bf2.3.~Remark.} {\sl Let\/ $R=\Bbb P_\Bbb CW$ be a geodesic, let\/ $w,w'\in W$ be a basis in an\/ $\Bbb R$-linear subspace\/
$W\subset V$ such that\/ $\langle W,W\rangle\subset\Bbb R$, and let\/ $h:=\langle-,w\rangle iw'-\langle-,w'\rangle iw$. Then\/ $R$ is
the real spine of the bisector\/ $B_h$.}

\medskip

{\bf Proof.} It is immediate that $h^*=h$, $\tr h=0$, $h$ has rank $2$, $R\subset B_h$, and $hW\subset iW$.

If $\ker h\cap W=0$, then $hW=iW$, implying that $R$ is the real spine of $B_h$.

If $\ker h\cap W\ne0$, then $R$ is euclidean, the focus $f\in W$ of $B_h$ is isotropic, and $B_h$ is parabolic. We~can assume that
$w'=w+f$. So, $h=\langle-,w\rangle if-\langle-,f\rangle iw$. Replacing $f,w\in W$ by $\Bbb R^*$-proportional ones, by Remark 2.2, we
can choose a $\Bbb C$-linear basis $p,hp,h^2p\in V$ such that $h^2p=f$, hence, $\langle p,f\rangle=1$. From
$hp=\langle p,w\rangle if-\langle p,f\rangle iw$, $\langle hp,p\rangle=0$, and $0\ne\langle p,f\rangle\in\Bbb R$, we conclude that
$\langle p,w\rangle\in\Bbb R$. Therefore, the~$\Bbb R$-span of $hp,ih^2p$ coincides with that of $iw,if$ because
$0\ne\langle w,w\rangle\in\Bbb R$
$_\blacksquare$

\medskip

Distinct projective lines are said to be {\it orthogonal\/} if their polar points are orthogonal.

\medskip

{\bf2.4.~Remark.} {\sl Two bisectors coincide if they have two common nonorthogonal slices.}

\medskip

{\bf Proof.} Such bisectors have the same real spines
$_\blacksquare$

\medskip

{\bf2.5.~Remark.} {\sl Let\/ $S$ be a projective line with the nonisotropic polar point\/ $p\notin\S$ and let\/ $f\in S$. Then there
exists a unique projective line\/ $S'$ orthogonal to\/ $S$ such that\/ $f\in S'$.}

\medskip

{\bf Proof.} We have $\S\not\ni p\in C:=\Bbb P_\Bbb Cf^\perp$. The point $p'$ polar to $S'$ is orthogonal to $p$ and belongs to $C$
$_\blacksquare$

\medskip

{\bf2.6.~Remark.} {\sl Let\/ $S$ be a slice of signature\/ $++$ or\/ $+-$ of a bisector\/ $B_h$. Then the projective line\/ $S'$
orthogonal to\/ $S$ and containing the focus\/ $f$ of\/ $B_h$ is a slice of\/ $B_h$.}

\medskip

{\bf Proof.} The point $p$ polar to $S$ is nonisotropic and belongs to the real spine $R$ of $B_h$. Denote by
$C:=\Bbb P_\Bbb Cf^\perp$ the complex spine of $B_h$ and by $p'$, the point polar to $S'$. As $p'\in C$ is orthogonal to
$p\in C\setminus\S$, the geodesic $R$ subject to $p\in R\subset C$ contains $p'$
$_\blacksquare$

\medskip

{\bf Proof of Remark 1.1.2.} If the bisectors are hyperbolic/spherical, the fact follows immediately from the considerations in 2.1.
So, we deal with the parabolic bisectors. The bisectors have the same focus and real spine. Therefore, by Remark 2.2, we can assume
that $\Bbb Rh_1p_1+\Bbb Rih_1^2p_1=\Bbb Rh_2p_2+\Bbb Rih_2^2p_2$, where the Gram matrices of $p_1,h_1p_1,h_1^2p_1$ and of
$p_2,h_2p_2,h_2^2p_2$ equal $G:=\left[\smallmatrix0&0&1\\0&1&0\\1&0&0\endsmallmatrix\right]$. In particular,
$\Bbb Rih_1^2p_1=\Bbb Rih_2^2p_2$. Replacing $h_2$ by an $\Bbb R^*$-proportional one and, if necessary, $p_2$ by $-p_2$, we obtain
$h_1^2p_1=h_2^2p_2$. Also, we~have $\pm h_1p_1=h_2p_2+2irh_2^2p_2$ for some $r\in\Bbb R$. Replacing $p_2$ by
$p'_2:=p_2+2irh_2p_2+2r^2h_2^2p_2$, we keep the same Gram matrix $G$ of $p'_2,hp'_2,h^2p'_2$ and obtain $h_1^2p_1=h_2^2p'_2$ and
$\pm h_1p_1=h_2p'_2$. Replacing $h_1$ by $\pm h_1$ and writing $p_2$ in place of $p'_2$, we get $h_1^2p_1=h_2^2p_2$ and
$h_1p_1=h_2p_2$. In the basis $p_1,h_1p_1,h_1^2p_1$, the~matrix of $h_2$ has the form
$h_2=\left[\smallmatrix0&0&0\\w&0&0\\r&1&0\endsmallmatrix\right]$. From the equality $h_2^tG=G\overline h_2$ (which is nothing but
$h_2^*=h_2$), we conclude that $w=1$ and $r\in\Bbb R$. Let $p_2=\left[\smallmatrix z\\y\\x\endsmallmatrix\right]$. From
$h_1p_1=h_2p_2$, we obtain $z=1$ and $y=-r$. So, $h_2=\left[\smallmatrix0&0&0\\1&0&0\\r&1&0\endsmallmatrix\right]$,
$p_2=\left[\smallmatrix1\\-r\\x\endsmallmatrix\right]$, and $h_2p_2=\left[\smallmatrix0\\1\\0\endsmallmatrix\right]$. It follows from
$\langle p_2,h_2p_2\rangle=0$ that $r=0$
$_\blacksquare$

\medskip

{\bf2.7.~Corollary.} {\sl Let\/ $B_{h_1},B_{h_2}$ be distinct parabolic bisectors with a common focus\/ $f\in\S$ and let\/
$R_1,R_2\subset C:=\Bbb P_\Bbb Cf^\perp$ denote their real spines. The following conditions are equivalent\/{\rm:}

\smallskip

$\bullet$ $B_{h_1}\cap B_{h_2}\cap\B=\varnothing${\rm;}

$\bullet$ $R_1,R_2\subset C$ are `parallel,' i.e., $R_1\cap R_2=\{f\}${\rm;}

$\bullet$ the\/ $\Bbb R$-span of\/ $h_1,h_2$ contains\/ $\langle-,f\rangle f$.}

\medskip

{\bf Proof.} The first two conditions are equivalent because every slice of $B_{h_i}$ except $C$ has signature $+-$ and the
intersection $B_{h_1}\cap B_{h_2}$ is the union of the projective lines whose polar points are in $R_1\cap R_2$.

Suppose that the $\Bbb R$-span of $h_1,h_2$ contains $\langle-,f\rangle f$. We can assume that $h_1+h_2=\langle-,f\rangle f$.
If~$p\in B_{h_1}\cap B_{h_2}\cap\B$, then $\big\langle\langle p,f\rangle f,p\big\rangle=0$, i.e., $\langle p,f\rangle=0$, a
contradiction.

Let $B_h$ be a bisector with focus $f$ and real spine $R\subset C$ and let $h_r:=h+r\langle-,f\rangle f$, $r\in\Bbb R$. If~$h_r$ has
rank $1$, then $h_r=\langle-,q\rangle p\ne0$. From $h_r^*=h_r$, $\tr h_r=0$, and $h_rf=0$, we conclude that
$\langle-,q\rangle p=\langle-,p\rangle q$, $\langle p,q\rangle=0$, and $\langle f,q\rangle p=0$. In other words, $p,q\in\Bbb Rf$,
hence, $h\in\Bbb R\langle-,f\rangle f$ has rank $\le1$, a contradiction. Therefore, $B_{h_r}$ is a bisector with focus $f$. For any
$p\in C$ different from $f$, we have $\langle h_rp,p\rangle=\langle hp,p\rangle+r\langle p,f\rangle\langle f,p\rangle=0$ for a
suitable $r\in\Bbb R$ because $\langle hp,p\rangle\in\Bbb R$ and $\langle p,f\rangle\ne0$. Consequently, the real spine $R_r$ of
$B_{h_r}$, a geodesic already known to be `parallel' to $R$, passes through~$p$. Thus, $R_r$, $r\in\Bbb R$, lists all geodesics in
$C$ that are `parallel' to $R$ (the geodesics on $C\setminus\{f\}$ can be seen as the geodesics on the usual euclidean plane)
$_\blacksquare$

\medskip

{\bf2.8.~Slice and meridional decompositions.} The {\it slice decomposition\/} of a bisector was already introduced in Definition
1.1.1: every point $p\in B_h$ different from the focus $f$ belongs to a unique slice $S_p\subset B_h$ and two different slices
intersect only in $f$.

Any projective line $L$ lying in $B_h$ is necessarily a slice. Otherwise, $f\notin L$ and $B_h$ contains the projective cone with apex
$f$ and base $L$, implying that $B_h=\Bbb P_\Bbb CV$, a contradiction.

As the real spine $R$ of $B_h$ is a geodesic, there is a $2$-dimensional $\Bbb R$-linear subspace $W\subset V$ such that
$R=\Bbb P_\Bbb CW$ and $0\ne\langle W,W\rangle\subset\Bbb R$.

Given a point $q\in B_h$ such that $f\ne q\notin C$, the point $hq\in R$ is the unique point in the real spine orthogonal to $q$ and
$hq\ne f$. We claim that $R$ and $q$ span the $\Bbb R$-plane $P_q=\Bbb P_\Bbb CW_q\subset B_h$, where $W_q$ stands for the
$\Bbb R$-span of $W$ and a representative $q\in V$ such that $\langle W,q\rangle\subset\Bbb R$. Indeed, if $B_h$ is
hyperbolic/spherical, then $W_q$ is the $\Bbb R$-span of $W\subset f^\perp$ and a representative $f\in V$.

So, we assume $B_h$ parabolic. In terms of the $\Bbb C$-linear basis $p,hp,h^2p\in V$ introduced in Remark~2.2, $W$ is the
$\Bbb R$-span of $hp,ih^2p$. Since $0\ne\langle ih^2p,q\rangle\in\Bbb R$, we can take $q=ip+rhp+ch^2p$, where $r,c\in\Bbb C$.
It~follows from $q\in B_h$ that $r\in\Bbb R$. Denoting $r':=\Re c$, we see that $W_q$ is the $\Bbb R$-span of $ip+r'h^2p,hp,ih^2p$. As
the Gram matrix of $ip+r'h^2p,hp,ih^2p$ equals $\left[\smallmatrix0&0&1\\0&1&0\\1&0&0\endsmallmatrix\right]$, we conclude that $P_q$
is an $\Bbb R$-plane and that $P_q\subset B_h$.

The $\Bbb R$-plane $P_q$ is spanned by $R\subset P_q$ and any $p\in P_q\setminus R$ different from $f$ because
$\langle p,R\rangle\ne0$. We~arrive at the {\it meridional decomposition\/} of $B_h$ : every point $q\in B_h$ different from the focus
$f$ and not lying on the complex spine $C$ belongs to a unique $\Bbb R$-plane $P_q\subset B_h$, called a {\it meridian\/} of $B_h$,
and~two different meridians intersect only in $R\cup\{f\}$, where $R$ stands for the real spine of $B_h$. Note that the union of all
meridians of a parabolic bisector $B_h$ does not contain the complex spine $C\subset B_h$. (The~other bisectors are unions of their
meridians.)

The intersection of a slice $S$ and a meridian $M$ is a geodesic passing through the focus because
$S\cap M=\Bbb P_\Bbb Cp^\perp\cap M$ is a geodesic, where $p\in R\subset M$ stands for the point polar to $S$.

\medskip

{\bf2.9.~Lemma.} {\sl Any geodesic contained in a bisector is contained in a slice or in a meridian of the bisector.}

\medskip

{\bf Proof.} Let $G\subset B_h$ be a geodesic in a bisector $B_h$.

Suppose that $B_h$ is hyperbolic/spherical. If $G$ intersects the real spine $R$ of $B_h$, then $G$ lies in a meridian of $B_h$. If
two distinct points of $G$ are in a same slice, then $G$ is contained in this slice. So,~we~pick distinct points $af+bw_1,f+w_2\in G$
such that $0\ne\langle af+bw_1,f+w_2\rangle\in\Bbb R$, the Gram matrix $\left[\smallmatrix r_1&r\\r&r_2\endsmallmatrix\right]$ of
$w_1,w_2\in R$ is real with $rr_2(r_1r_2-r^2)\ne0$, and $0\ne a,b\in\Bbb C$, where $f$ stands for the focus of~$B_h$.

For any $x\in\Bbb R$, we have $af+bw_1+x(f+w_2)\in B_h$. Therefore, $bw_1+xw_2\in R$ for infinitely many $x\in\Bbb R$. Since
$\langle bw_1+xw_2,w_2\rangle=rb+r_2x$, we obtain
$r_2\Big\langle\displaystyle\frac{bw_1+x_1w_2}{rb+r_2x_1},\frac{bw_1+x_2w_2}{rb+r_2x_2}\Big\rangle\in\Bbb R$ for distinct
$x_1,x_2\in\Bbb R$. This means that $1+\displaystyle\frac{(r_1r_2-r^2)|b|^2}{r^2|b|^2+rr_2(x_2b+x_1\overline b)+r_2^2x_1x_2}\in\Bbb R$
and implies $\Bbb R\ni x_2b+x_1\overline b=(x_2-x_1)b+2x_1\Re b$, i.e., $b\in\Bbb R$. It follows from
$\langle af+bw_1,f+w_2\rangle\in\Bbb R$ that $a\langle f,f\rangle+br\in\Bbb R$ and, hence, $a\in\Bbb R$. As $a,b\in\Bbb R$, the points
$af+bw_1,f+w_2$ lie in a meridian of $B_h$. Consequently, $G$ is contained in this meridian.

Suppose that $B_h$ is parabolic. Let $p,hp,h^2p\in V$ be a $\Bbb C$-linear basis with the Gram matrix
$\left[\smallmatrix0&0&1\\0&1&0\\1&0&0\endsmallmatrix\right]$ (as~the one introduced in Remark 2.2). Assuming that $G$ does not lie in
the complex spine $C$ (spanned by $hp,h^2p$), we~pick distinct nonorthogonal points $q_1,q_2\in G\setminus C$. We can assume that
$q_1=uip+ua_1hp+b_1h^2p$, $q_2=ip+a_2hp+b_2h^2p$ and $0\ne\langle q_1,q_2\rangle\in\Bbb R$, where $a_i,b_i,u\in\Bbb C$ and $|u|=1$. It
follows from $\langle hq_i,q_i\rangle=0$ that $a_i\in\Bbb R$. Since $\langle q_1,q_2\rangle\in\Bbb R$, we obtain
$u(a_1a_2+i\overline b_2)-ib_1\in\Bbb R$.

For any $x\in\Bbb R$, we have $\big\langle h(q_1+xq_2),q_1+xq_2\big\rangle=0$, i.e., $(a_1-a_2)x\Im u=0$. So, either $a_1=a_2$ or
$u=\pm1$.

In the case of $u=\pm1$, we can assume that $u=1$ and derive $\Re b_1=\Re b_2$ from $u(a_1a_2+i\overline b_2)-ib_1\in\Bbb R$. It
follows that $q_1-q_2=(a_1-a_2)hp+(b_1-b_2)h^2p\in\Bbb Rhp+\Bbb Rih^2p$. In other words, $q_1-q_2$ belongs to the real spine of $B_h$,
implying that $q_1,q_2$ are in a meridian of $B_h$ and that $G$ lies in this meridian.

In the case of $a_1=a_2$, we obtain $\overline uq_1-q_2\in\Bbb Ch^2p$. As $h^2p$ is the focus of $B_h$, the points $q_1,q_2$ are in a
same slice of $B_h$. Therefore, $G$ lies in this slice
$_\blacksquare$

\medskip

{\bf2.10.~Corollary.} {\sl Any\/ $\Bbb R$-plane contained in a bisector is a meridian of the bisector.}

\medskip

{\bf Proof.} Let $P\subset B_h$ be an $\Bbb R$-plane. It suffices to observe that $P$ contains the real spine $R$ of $B_h$, i.e., that
the intersection $P\cap R$ has at least $3$ points. Let $p\in P\setminus R$ and let $p\in G_i\subset P$, $i=1,2,3$, be distinct
geodesics. By Lemma 2.9, $G_i$ is contained in a meridian $R_i$. Since $R\subset R_i$, we obtain $G_i\cap R\ne\varnothing$
$_\blacksquare$

\medskip

{\bf2.11.~Normal vector.} Let $p\in B_h$ be a nonisotropic point in the bisector $B_h$ different from the focus $f$ of $B_h$,
$p\ne f$. Then $n_p:=\langle-,p\rangle hp$ is a normal vector to $B_h$ at $p$. Indeed, $n_p\ne0$ because $p\ne f$. Since $p\in B_h$,
we have $\langle hp,p\rangle=0$. Let $t:=\langle-,p\rangle v\in\T_p\Bbb P_\Bbb CV$, $\langle v,p\rangle=0$, be a tangent vector to
$B_h$ at $p$. Taking the derivative $\frac d{d\varepsilon}\Big|_{\varepsilon=0}$ of
$\big\langle h(p+\varepsilon tp),p+\varepsilon tp\big\rangle=0$, we obtain $\Re\langle hp,v\rangle=0$. In other words, $n_p$ is
orthogonal to any vector tangent to $B_h$ at $p$. In particular, the bisector $B_h$ is smooth at~its (nonisotropic) points different
from the focus $f$.

\medskip

{\bf2.12.~Lemma.} {\sl Let\/ $B_h$ be a bisector. Denote by\/ $C$ and\/ $R$ the complex and real spines of\/ $B_h$. Then, for any\/
$c_1\in C\setminus R$, there exists a unique\/ $c_2\in C\setminus R$, called the\/ {\rm reflection} of\/ $c_1$ in\/ $R$, such that\/
$h=\langle-,c_1\rangle\frac{c_1}{r_1}+\langle-,c_2\rangle\frac{c_2}{r_2}$ for suitable\/ $0\ne r_1,r_2\in\Bbb R$ depending on
representatives\/ $c_1,c_2\in V$. If\/~$c_1=hq_1$, then\/ $r_1:=\langle q_1,hq_1\rangle$ and\/ $r_2:=\langle q_2,hq_2\rangle$, where\/
$c_2:=hq_2\in C$ is uniquely determined by\/ $\langle q_1,c_2\rangle=0$. Moreover, $c_1,c_2$ have a same signature and\/ $r_1r_2<0$.
Therefore, $\pm h=\langle-,c_1\rangle c_1-\langle-,c_2\rangle c_2$ for suitable representatives\/ $c_1,c_2\in V$.}

\medskip

{\bf Proof.} Any $c\in C$ has the form $c=hp$, where $p$ is different from the focus $f$ of $B_h$. Since $hp$ is the point polar to
the slice $S_p$ of $p\in B_h\setminus\{f\}$, we conclude that $c\notin R$ iff $p\notin B_h$.

Let $c_1=hp_1$. As $p_1\ne f$ and $C=\Bbb P_\Bbb Cf^\perp$, we can define $c$ as $\{c\}:=\Bbb P_\Bbb Cp_1^\perp\cap C$. For some
$p\ne f$, we~have $c=hp$. It follows from $\langle hp,p_1\rangle=0$ that $\langle p,c_1\rangle=\langle p,hp_1\rangle=0$.

Suppose that $c\in R$. Then $\langle p,hp\rangle=0$ and $hp\ne hp_1$ as, otherwise,
$\langle p_1,hp_1\rangle\in\langle p_1,\Bbb Chp\rangle=0$, contradicting $c_1\notin R$. So, $hV$ is the $\Bbb C$-span of $hp,hp_1$.
From $\langle p,hp\rangle=\langle p,hp_1\rangle=0$, we deduce that $p=f$, a contradiction. Therefore, $c\in C\setminus R$.

For $h$ of the form indicated in Lemma 2.12, we obtain $c=hp=\langle p,c_2\rangle\frac{c_2}{r_2}$. Thus, $c_2$ is unique.

It follows from $\langle p,hp_1\rangle=\langle p_1,hp\rangle=0$ and $\langle f,hV\rangle=0$ that
$h':=\frac{\langle-,hp_1\rangle hp_1}{\langle p_1,hp_1\rangle}+\frac{\langle-,hp\rangle hp}{\langle p,hp\rangle}$ and $h$ coincide on
$p_1,p,f$. Since $p_1,p,f$ are not in a same projective line, $h'=h$.

If $f\in\B\cup\S$, then $c_1,c_2$ are positive and $\tr h=0$ implies $r_1r_2<0$.

Consider the remaining case of a hyperbolic $B_h$. Pick a basis $v_1,f,v_2\in V$ of eigenvectors of $h$ corresponding to the
eigenvalues $-ri,0,ri$, $0\ne r\in\Bbb R$, with the Gram matrix $\left[\smallmatrix0&0&1\\0&1&0\\1&0&0\endsmallmatrix\right]$. Any
$c\in C\setminus R$ has the form $c=h(r^{-1}iav_1-r^{-1}iv_2)=av_1+v_2$ with $a\in\Bbb C\setminus\Bbb R$. Let $c_1=av_1+v_2$. Then
$c_2=\overline av_1+v_2$ because $\langle r^{-1}iav_1-r^{-1}iv_2,\overline av_1+v_2\rangle=0$. As
$\langle c_1,c_1\rangle=2\Re a=\langle c_2,c_2\rangle$, the points $c_1,c_2$ have the same signature. The matrices of
$\langle-,c_1\rangle c_1$, $\langle-,c_2\rangle c_2$, $h$ equal
$\left[\smallmatrix a&0&|a|^2\\0&0&0\\1&0&\overline a\endsmallmatrix\right]$,
$\left[\smallmatrix\overline a&0&|a|^2\\0&0&0\\1&0&a\endsmallmatrix\right]$,
$\left[\smallmatrix-ri&0&0\\0&0&0\\0&0&ri\endsmallmatrix\right]$, implying the rest
$_\blacksquare$

\bigskip

\centerline{\bf3.~Proof of Proposition 1.2.3}

\medskip

This section contains a series of elementary technical linear algebra facts required mostly in the proof of Proposition 1.2.3.

\medskip

{\bf3.1.~Lemma.} {\sl Let\/ $\Bbb K$ be a field of characteristic\/ $\ne2$ and let\/ $K,W,A$ be\/ $\Bbb K$-linear spaces of
dimensions\/~$3,2,3$, respectively, such that\/ $K\subset W\otimes_\Bbb KA$ and\/ $K\cap(w\otimes A)\ne0$ for every\/ $0\ne w\in W$.
Then\/ $W\otimes a\subset K$ for some\/ $0\ne a\in A$ or\/ $K\subset W\otimes_\Bbb KA'$ for some\/ $\Bbb K$-linear subspace\/
$A'\subset A$ with\/ $\dim_\Bbb KA'\le2$.}

\medskip

{\bf Proof.} Let $w_1,w_2\in W$ be a $\Bbb K$-linear basis. For some $0\ne a_1,a_2,a,a'\in A$, we have
$$w_1\otimes a_1,w_2\otimes a_2,(w_1+w_2)\otimes a,(w_1-w_2)\otimes a'\in K.$$

If $w_1\otimes a_1,w_2\otimes a_2,(w_1+w_2)\otimes a$ are linearly dependent, we have
$w_1\otimes c_1a_1+w_2\otimes c_2a_2=(w_1+w_2)\otimes a$ for some $c_1,c_2\in\Bbb K$, implying $a=c_1a_1=c_2a_2$ and $c_1,c_2\ne0$.
So, $W\otimes a\subset K$ and we are done.

Therefore, we assume that $w_1\otimes a_1,w_2\otimes a_2,(w_1+w_2)\otimes a$ form a basis in $K$. It follows that
$(w_1-w_2)\otimes a'=w_1\otimes c_1a_1+w_2\otimes c_2a_2+(w_1+w_2)\otimes ca$ for suitable $c_1,c_2,c\in\Bbb K$. Hence,
$a'=c_1a_1+ca$ and $-a'=c_2a_2+ca$, implying $c_1a_1+c_2a_2+2ca=0$. If $a_1,a_2,a$ are linearly independent, we obtain $c_1=c_2=c=0$
and $a'=0$, a contradiction. Consequently, $a_1,a_2,a\in A'\subset A$ and $K\subset W\otimes A'$, where $\dim_\Bbb KA'\le2$
$_\blacksquare$

\medskip

{\bf3.2.~Corollary.} {\sl Let\/ $W_0\subset\Lin_\Bbb C(V,V)$ be a\/ $3$-dimensional\/ $\Bbb R$-linear subspace such that\/
$\det W_0=0$. Then there exist\/ $\Bbb C$-linear subspaces\/ $V_1,V_2\subset V$ such that\/ $\dim_\Bbb CV_i=i$ for\/ $i=1,2$ and
either\/ $W_0V\subset V_2$ or\/ $W_0V_1=0$ or\/ $W_0V_2\subset V_1$.}

\medskip

{\bf Proof.} Take $\Bbb K:=\Bbb C$, $W:=\Bbb CW_0$, and $A:=V$ in Lemma 3.1. The rule $\varphi:h\otimes a\mapsto ha$ defines a
$\Bbb C$-linear map $\varphi:W\otimes_\Bbb CA\to V$. Denote $K:=\ker\varphi$. If the image of $\varphi$ lies in a $2$-dimensional
$\Bbb C$-linear subspace $V_2\subset V$, we have $WV\subset V_2$, hence, $W_0V\subset V_2$. If $\dim_\Bbb CW=1$, we can take
$V_1:=\Bbb Cf$ for $0\ne f\in\ker h_0$ and $0\ne h_0\in W_0$, thus getting $W_0V_1\subset WV_1=0$. So, we can assume that
$\dim_\Bbb CW=2$ and $\dim_\Bbb CK=3$. In order to verify the conditions of Lemma 3.1, it suffices to observe that $\det W=0$.
Let $h_1,h_2$ be an $\Bbb R$-linear basis of $W_0$. Then the polynomial $p(x_1,x_2):=\det(h_1x_1+h_2x_2)$ vanishes for all
$x_1,x_2\in\Bbb R$. So, it vanishes identically. By Lemma 3.1, $Wf=0$ for some $0\ne f\in V$ or $K\subset W\otimes V_2$ for a suitable
$2$-dimensional $\Bbb C$-linear subspace $V_2\subset V$. In the first case, we have $W_0V_1=0$, where $V_1:=\Bbb Cf$. In the second
case, $\dim_\Bbb CWV_2=4-3=1$, hence, $W_0V_2\subset V_1$, where $V_1:=WV_2$
$_\blacksquare$

\medskip

{\bf3.3.~Lemma.} {\sl Suppose that any nonnull\/ $\Bbb C$-linear map from an $\Bbb R$-linear subspace\/ $W\subset\Lin_\Bbb C(V,V)$ has
rank\/ $2$. Then either\/ $\dim_\Bbb CWV\le2$ or\/ $Wf_0=0$ for some\/ $0\ne f_0\in V$ or there exists a unique up to\/
$\Bbb C^*$-proportionality\/ $\Bbb R$-linear embedding\/ $f:W\hookrightarrow V$ such that\/ $h(fh)=0$ for all\/ $h\in W$ and, hence,
$h(fh')+h'(fh)=0$ for all\/ $h,h'\in W$.}

\medskip

{\bf Proof.} It is immediate that the identity $h(fh)=0$, $h\in W$, for an $\Bbb R$-linear map $f:W\to V$ implies the identity
$h(fh')+h'(fh)=0$, $h,h'\in W$. (One may substitute $h+h'$ in $h(fh)=0$.)

Let $0\ne h_0\in W$ and $0\ne f_0\in V$ be such that $h_0f_0=0$. Define $\Bbb R$-linear subspaces
$$N(h_0):=\{h\in W\mid\ker h_0\subset\ker h\},\qquad N'(h_0):=\{h\in W\mid hV\subset h_0V\}.$$
If $N(h_0)=W$ or $N'(h_0)=W$, we are done because $\dim_\Bbb Ch_0V=2$.

So, we assume that $U:=W\setminus\big(N(h_0)\cup N'(h_0)\big)\ne\varnothing$. Let $u\in U$. We can apply Corollary 3.2 to the
$\Bbb R$-span $W_0$ of $h_0,u$ because $h_0\in N(h_0)$ implies $\dim_\Bbb RW_0=2$. By the choice of $u$, there exist $\Bbb C$-linear
subspaces $V_i\subset V$ such that $\dim_\Bbb CV_i=i$ for $i=1,2$ and $W_0V_2\subset V_1$. Since each of $h_0,u$ has rank $2$,
we~obtain $\ker h_0,\ker u\subset V_2$. From $u\notin N(h_0)$, we conclude that $\ker h_0\ne\ker u$, hence, $V_2=\ker h_0+\ker u$ and
$V_1=h_0(\ker u)=u(\ker h_0)$. In other words, there exists a unique $0\ne\psi u\in\ker u$ such that $h_0(\psi u)+uf_0=0$. We have
defined a map $\psi:U\to V$ such that $\psi u\ne0$, $u(\psi u)=0$, and $h_0(\psi u)+uf_0=0$ for all $u\in U$. It follows from the
uniqueness of $\psi u$ satisfying the equality $h_0(\psi u)+uf_0=0$ that $\psi(ru)=r\psi u$ and $\psi(u+u')=\psi u+\psi u'$ for all
$0\ne r\in\Bbb R$ and $u,u'\in U$ such that $u+u'\in U$. One can readily observe that $\psi$ is continuous.

We pick an $\Bbb R$-linear basis $b_1,b_2,\dots,b_n\in U$ of $W$ and define an $\Bbb R$-linear map $f:W\to V$ by the rule
$f\big(\sum_{i=1}^nr_ib_i\big):=\sum_{i=1}^nr_i\psi b_i$ for all $r_1,r_2,\dots,r_n\in\Bbb R$. Then $f$ and $\psi$ coincide on a
subset $U_0\subset U$ dense in~$W$. For example, one can take
$$U_0:=\Big\{\sum_{i=1}^nr_ib_i\mid r_k\ne0\text{ and }\sum_{i=1}^kr_ib_i\notin N(h_0)\cup N'(h_0)\text{ for all }1\le k\le n\Big\}.$$
As $\psi$ is continuous, we obtain $\psi=f|_U$. Hence, $h(fh)=0$ and $fu\ne0$ for all $h\in W$ and $u\in U$.

We will show that $f$ is unique up to $\Bbb C^*$-proportionality. Let $f':W\to V$ be another nonnull $\Bbb R$-linear map such that
$h(f'h)=0$ for all $h\in W$. For some $h\in W$, we have $0\ne fh,f'h\in\ker h$. Replacing $f'$ by a $\Bbb C^*$-proportional one, we
obtain $fh=f'h\ne0$. The equalities $h(fh')+h'(fh)=0$, $h(f'h')+h'(f'h)=0$, and $fh=f'h$ imply $h(fh'-f'h')=0$. So, $fh'=f'h'$ if
$\ker h'\cap\ker h=0$. It remains to observe that the $h'\in W$ subject to $\ker h'\cap\ker h=0$ form a dense subset in $W$ unless
$W\ker h=0$.

In order to show that $f$ is an embedding, it suffices to construct $f$ that does not vanish on a given $0\ne h\in W$. So, it suffices
to find $0\ne h_0\in W$ such that $h\notin N(h_0)\cup N'(h_0)$. Suppose that $h\in N(h_0)\cup N'(h_0)$ for any $0\ne h_0\in W$. This
means that $\ker h_0=\ker h$ or $h_0V=hV$ for any $0\ne h_0\in W$. In other words, $h_0\in N(h)\cup N'(h)$ for any $h_0\in W$
$_\blacksquare$

\medskip

{\bf3.4.~Remark.} {\sl Suppose that\/ $h\in\Lin_\Bbb C(V,V)$ satisfies\/ $h^*=h$ and\/ $\tr h=0$. If\/ $h$ has rank\/ $\le1$, then\/
$h\in\Bbb R\langle-,f_0\rangle f_0$, where\/ $f_0\in V$ is isotropic.}

\medskip

{\bf Proof.} Suppose that $h$ has rank $1$. Then $h=\langle-,f_0\rangle v$ for some $f_0,v\in V$. From $h^*=h$ and $\tr h=0$, we
conclude that $\langle-,f_0\rangle v=\langle-,v\rangle f_0$ and $\langle v,f_0\rangle=0$. So, $v=rf_0$ with $f_0\in\S$ and
$0\ne r\in\Bbb R$
$_\blacksquare$

\medskip

The following corollary follows directly from Remark 3.4.

\medskip

{\bf3.5.~Corollary.} {\sl Let\/ $U\subset V$ be a $2$-dimensional\/ $\Bbb C$-linear subspace. Denote
$$W_U:=\big\{h\in{\Lin}_\Bbb C(V,V)\mid\tr h=0,\ h^*=h,\ hU\subset U^\perp\big\}.\leqno{\bold{(3.6)}}$$
If\/ $U$ is nondegenerate, then every\/ $0\ne h\in W_U$ has rank\/ $2$. If\/ $U$ is degenerate, then\/ $\Bbb R\langle-,f_0\rangle f_0$
is the set of all elements in\/ $W_U$ of rank\/ $\le1$, where\/ $\Bbb Cf_0:=U^\perp$}
$_\blacksquare$

\medskip

{\bf Proof of Proposition 1.2.3.} If $\dim_\Bbb CWV\le2$, then there exists $0\ne f_0\in V$ such that $\langle WV,f_0\rangle=0$. Since
$W^*=W$, we obtain $\langle V,Wf_0\rangle=0$, implying $Wf_0=0$, a contradiction. By Lemma 3.3, we get an $\Bbb R$-linear embedding
$f:W\hookrightarrow V$ such that $h(fh)=0$ for all $h\in W$.

As $Wf_0\ne0$ for all\/ $0\ne f_0\in V$, we have $\dim_\Bbb CU\ge2$.

Suppose that $\dim_\Bbb CU=2$. Let $0\ne h\in W$. Then $\dim_\Bbb ChU=1$ because $0\ne\ker h\subset U$. Let $h'\in W$. If
$fh'\notin\Bbb Cfh$, then $h(fh')\ne0$ and
$$\big\langle h(fh'),fh'\big\rangle=-\big\langle h'(fh),fh'\big\rangle=-\big\langle fh',h'(fh')\big\rangle=0,$$
i.e., $\langle hU,fh'\rangle=0$. If $fh'=c(fh)$ for $c\in\Bbb C$, then
$$\langle hU,fh'\rangle=\overline c\langle hU,fh\rangle=\overline c\big\langle U,h(fh)\big\rangle=0.$$
Consequently, $\langle hU,fW\rangle=0$ and $\langle hU,U\rangle=0$. This means that $h\in W_U$.

Suppose that $U=V$. Let $f_i:=fh_i$, $i=1,2,3$, be a $\Bbb C$-linear basis in $V$, $h_i\in W$. Denote $e_i:=h_{i+1}f_{i+2}$ (the
indices are modulo $3$). By Lemma 3.3, $h_if_{i+1}+h_{i+1}f_i=0$. Since $h_i^*=h_i$ and $h_if_i=0$, we obtain
$$\langle f_{i-1},e_i\rangle=\langle h_{i+1}f_{i+2},f_{i+2}\rangle=-\langle h_{i+2}f_{i+1},f_{i+2}\rangle=-\langle
f_{i+1},h_{i+2}f_{i+2}\rangle=0,$$
$$\langle f_{i+1},e_i\rangle=\langle h_{i+1}f_{i+1},f_{i+2}\rangle=0,\qquad\langle f_i,e_i\rangle=\langle
h_{i+1}f_i,f_{i+2}\rangle=-\langle h_if_{i+1},f_{i+2}\rangle=-\langle e_{i-1},f_{i-1}\rangle.$$
Denoting $c_i:=\langle f_i,e_i\rangle$, we get $c_i=-\overline c_{i-1}$ and, consequently, $c_i=c_{i-2}$. So,
$0\ne ri:=c_1=c_2=c_3\in\Bbb Ri$. Let $G:=[g_{ij}]$ denote the Gram matrix of the $f_i$'s and let $A:=[a_{ij}]$ be the matrix
expressing the $e_j$'s in terms of the basis $f_1,f_2,f_3$, i.e., $g_{ij}:=\langle f_i,f_j\rangle$ and $e_j=\sum_ia_{ij}f_i$. Then, as
we have shown, $A^tG=ri$. Since $h_if_i=0$, $h_if_{i+1}=e_{i-1}$, and $h_if_{i-1}=-h_{i-1}f_i=-e_{i-2}$, in the basis $f_1,f_2,f_3$,
we have
$$h_1=\left[\smallmatrix0&a_{13}&-a_{12}&\\0&a_{23}&-a_{22}\\0&a_{33}&-a_{32}\endsmallmatrix\right],\qquad
h_2=\left[\smallmatrix-a_{13}&0&a_{11}&\\-a_{23}&0&a_{21}\\-a_{33}&0&a_{31}\endsmallmatrix\right],\qquad
h_3=\left[\smallmatrix a_{12}&-a_{11}&0\\a_{22}&-a_{21}&0\\a_{32}&-a_{31}&0\endsmallmatrix\right].$$
It follows from $\tr h_i=0$ that $A$ is symmetric, $A^t=A$. We conclude that $iG$ is symmetric. Therefore, $G$ is real.

Denote by $W'$ the $\Bbb R$-span of $h_1,h_2,h_3$. We have shown that $fW'$ is a totally real subspace. Let~$h\in W\setminus W'$. We
can assume that $fh\in ifW'$. In other words, $fh=ifh'_1$ with $0\ne h'_1\in W'$. Pick an $\Bbb R$-linear basis $h'_1,h'_2,h'_3\in W'$
such that $\langle fh'_1,fh'_2\rangle\ne0$. We get a $\Bbb C$-linear basis $fh,fh'_2,fh'_3\in V$ because $fh'_1,fh'_2,fh'_3\in V$ is a
$\Bbb C$-linear basis. Applying the above to $h,h'_2,h'_3$ in place of $h_1,h_2,h_3$, we~conclude that
$\langle fh,fh'_2\rangle\in\Bbb R$. On the other hand, $\langle fh,fh'_2\rangle=i\langle fh'_1,fh'_2\rangle\in\Bbb Ri$, a
contradiction
$_\blacksquare$

\bigskip

\centerline{\bf4.~Linear families of bisectors and their intersections}

\medskip

In order to study intersections of finitely many bisectors, we need to analyze the most simple families of bisectors, the linear
families. For the vast majority of such families, we show that they are determined by their intersection (Proposition 1.2.4) and that
there is an $\Bbb R$-linear embedding that maps a bisector of the family into its focus (see Proposition 1.2.3). In addition to the
proofs of Remark 1.2.1 and Proposition 1.2.4, we describe the intersection of bisectors of a linear family and study the points where
it is transversal (Lemmas 4.4.1, 4.4.2, 4.4.3, 4.4.4, and 4.5.2). We also study when $2$ bisectors are transversal at a common
negative point (Lemma 4.6). The section ends with a version of the G.~Giraud rigidity theorem (Lemma 4.7 and Corollary 4.8). The
material of this section includes some facts from [Gol] presented in a more general form and under a different angle of view.

\medskip

{\bf Proof of Remark 1.2.1.} As the focus is a singular point of a bisector, Remark 1.2.1 holds when $p$ is a focus of some $B_{h_i}$,
i.e., if $h_ip=0$ for some $i$. So, we assume that $p$ is not a focus of $B_{h_i}$ for all $i$. By~2.11, $n_i:=\langle-,p\rangle h_ip$
is a normal vector to $B_{h_i}$ at $p$. As $p\notin\S$, the $n_i$'s are $\Bbb R$-linearly independent iff the $h_ip$ are
$\Bbb R$-linearly independent
$_\blacksquare$

\medskip

{\bf4.1.~Remark.} {\sl Let\/ $U\subset V$ be a\/ $2$-dimensional\/ $\Bbb C$-linear subspace. Then\/ $W_U$ given by\/ {\rm(3.6)} equals
$$W_U=\big\{h\in{\Lin}_\Bbb C(V,V)\mid\tr h=0,\ h^*=h,\ \langle hs,s\rangle=0\text{ \rm for all }s\in S\big\},\leqno{\bold{(4.2)}}$$
where\/ $S:=\Bbb P_\Bbb CU$.}

\medskip

{\bf Proof.} The $\Bbb R$-linear subspace $W_U$ defined by (4.2) obviously contains the one defined by (3.6).
If~$\langle hu,u\rangle=0$ for all $u\in U$, then $\big\langle h(u+u'),u+u'\big\rangle=0$ for all $u,u'\in U$ and
$\Re\langle hu,u'\rangle=0$, implying $\langle hu,U\rangle=0$ because $U$ is a $\Bbb C$-linear subspace
$_\blacksquare$

\medskip

Using Proposition 1.2.3, we are going to describe all linear families of bisectors and their bases.

\medskip

{\bf4.3.~Confocal linear families of bisectors.} A linear family $W$ of bisectors is {\it confocal\/} with the common focus
$f_0\in\Bbb P_\Bbb CV$ if $Wf_0=0$. A $2$-dimensional confocal linear family of bisectors is called a {\it confocal line\/} of
bisectors.

\medskip

{\bf4.3.1.~Lemma.} {\sl Let\/ $S$ be a projective line with the polar point\/ $p\not\in\S$ and let\/ $S'$ be a projective line
orthogonal to\/ $S$. Then
$$W:=\big\{h\in{\Lin}_\Bbb C(V,V)\mid\tr h=0,\ h^*=h,\ \langle hs,s\rangle=0\text{ \rm for all }s\in S\cup
S'\big\}.\leqno{\bold{(4.3.2)}}$$
is a confocal line of bisectors with the common focus\/ $f_0$ and the base\/ $B_W=S\cup S'$, where $S\cap S'=\{f_0\}$. Moreover,
$0\ne h\in W$ iff\/ $f_0$ is the focus of the bisector\/ $B_h$ and\/ $p$ belongs to the real spine of\/ $B_h$.

Every confocal line of bisectors has the above form.}

\medskip

{\bf Proof.} Let $W$ be given by (4.3.2). By Corollary 3.5 and Remark 4.1, every $0\ne h\in W$ has rank $2$. Furthermore, $p$ belongs
to the real spine of $B_h$ and, by 2.8, $S,S'$ are slices of $B_h$, implying that $f_0$ is the focus of $B_h$. As $S\cup\{p\}$
contains a negative point and $S\cup\{p\}\subset S\cup S'$, we conclude that $W$ is a linear family of bisectors.

Conversely, let $B_h$ be a bisector with the focus $f_0$ whose real spine contains $p$. Then $S\subset B_h$. By~Remark 2.6, $S'$ is a
slice of $B_h$.

Denote by $C:=\Bbb P_\Bbb Cf_0^\perp$ the common complex spine of the bisectors of the family. The point $p'\in C$ polar to $S'$ is a
unique point in $C$ orthogonal to $p\in C\setminus\S$. Clearly, $p'\ne p$.

Every $B_h$, $0\ne h\in W$, is completely determined by its real spine $R$ and $p\in R\subset C$. The geodesics $G$ subject to the
condition $p\in G\subset C$ form a $1$-parameter family and the intersection of all such geodesics equals $\{p,p'\}$. This implies
$\dim_\Bbb RW=2$ and $B_W=S\cup S'$.

\smallskip

Let $h_1,h_2\in W$ be an $\Bbb R$-linear basis of a confocal line $W$ of bisectors with the common focus~$f_0$. By~Remarks 2.5 and
2.6, it suffices to show that the bisectors $B_{h_1},B_{h_2}$ have a common slice of signature~$+-$. This is clear if $f_0\notin\B$
because $B_W=B_{h_1}\cap B_{h_2}$ and $B_W\cap\B\ne\varnothing$. Suppose that $f_0\in\B$. The real spines $R_1,R_2\subset C$ of
$B_{h_1},B_{h_2}$ intersect because any two geodesics lying in the projective line $C:=\Bbb P_\Bbb Cf_0^\perp$ (the common complex
spine of the bisectors) of signature $++$ do so
$_\blacksquare$

\medskip

{\bf4.3.3.~Remark.} {\sl Let\/ $W$ be a confocal linear family of bisectors with a nonnegative common focus\/~$f_0$, $f_0\notin\B$.
If\/ $\dim_\Bbb RW\ge2$, then\/ $W$ is a confocal line.}

\medskip

{\bf Proof.} Pick a point $p\in\B\cap B_W$ and denote by $S$ the projective line of signature $+-$ spanned by $f_0,p$. By Remark 2.5,
there is a unique projective line $S'$ orthogonal to $S$ such that $f_0\in S'$. By Remark 2.6, any bisector with the focus $f_0$ that
contains $S$ necessarily contains $S'$. By Lemma 4.3.1, all such bisectors form a confocal line. In other words,
$W':=\big\{h\in\Lin_\Bbb C(V,V)\mid\tr h=0,\ h^*=h,\ hf_0=0,\ \langle hp,p\rangle=0\big\}$ is a confocal line of bisectors. It remains
to observe that $W\subset W'$
$_\blacksquare$

\medskip

{\bf4.3.4.~Remark.} {\sl Let\/ $W$ be a confocal linear family of bisectors with negative common focus\/ $f_0$, $f_0\in\B$. Then\/
$\dim_\Bbb RW\le4$. If\/ $\dim_\Bbb RW\ge3$, then\/ $B_W=\{f_0\}$.}

\medskip

{\bf Proof.} Clearly, $\hat W:=\big\{h\in\Lin_\Bbb C(V,V)\mid\tr h=0,\ h^*=h,\ hf_0=0\big\}$ is a maximal confocal linear family of
bisectors. As $i\hat W$ is the Lie algebra of the stabilizer $\Stab f_0\simeq\U_2\Bbb C$ of $f_0$ in $\SU V$, we conclude that
$\dim_\Bbb R\hat W=4$.

By Lemma 4.3.1, the base $B_{W'}$ of any confocal line $W'\subset\hat W$ is the union of orthogonal projective lines $S,S'$ of
signature $+-$, $B_{W'}=S\cup S'$, and $S$ completely determines $W'$ ($f_0\in\B$ is fixed).

Let $W\subset\hat W$ be a linear family of dimension $\ge3$. Then $B_W=\{f_0\}$ because, otherwise, $B_W$ would contain a projective
line $S$ of signature $+-$, implying that $W$ is included in the confocal line $W'\subset W$ determined by $S$
$_\blacksquare$

\medskip

{\bf4.4.~Nonconfocal linear families of bisectors.} Consider a linear family $W$ of bisectors having no common focus. By Proposition
1.2.3, we have an $\Bbb R$-linear embedding $f:W\hookrightarrow V$ such that $h(fh)=0$ for all $h\in W$. It follows from
$$\big\langle h(fh'),fh'\big\rangle=-\big\langle h'(fh),fh'\big\rangle=-\big\langle fh,h'(fh')\big\rangle=0,$$
$h,h'\in W$, that $\Bbb P_\Bbb CfW\subset B_W$.

\medskip

{\bf4.4.1.~Lemma.} {\sl Let\/ $W$ be a linear family of bisectors whose foci do not lie in a same projective line. Then the foci of
the bisectors of the family form an\/ $\Bbb R$-plane\/ $P$, $\dim_\Bbb RW=3$, $B_W=P$, and
$$W=\big\{h\in{\Lin}_\Bbb C(V,V)\mid\tr h=0,\ h^*=h,\ \langle hp,p\rangle=0\text{ \rm for all }p\in P\big\}.$$
Conversely, given an\/ $\Bbb R$-plane\/ $P$, the above formula defines a nonconfocal linear family of bisectors whose foci
constitute\/ $P$.

No\/ $3$ bisectors of such a family are transversal at any\/ $p\in P\setminus\S$.}

\medskip

{\bf Proof.} By Proposition 1.2.3 and the beginning of 4.4, $\Bbb CfW=V$, $\dim_\Bbb RW=3$, and the bisectors~$B_h$, $0\ne h\in W$,
have a common $\Bbb R$-plane $P:=\Bbb P_\Bbb CfW$. Moreover, $P$ determines $W$ and all bisectors that contain $P$ have the form
$B_h$, $0\ne h\in W$. Indeed, by 2.1, every bisector $B$ is determined by its real spine $R$ which can be an arbitrary geodesic. If
$P\subset B$, then, by Corollary 2.10, $R\cup\{f\}\subset P$, where $f$ stands for the focus of $B$. As
$R=\Bbb P_\Bbb Cf^\perp\cap P$, a bisector $B$ containing $P$ is completely determined by its focus $f\in P$ which can be an arbitrary
point in $P$. We call such a family $W$ an {\it$\Bbb R$-plane\/} of bisectors.

Let $p\notin P$ and let $\Gamma:=P\cap\Bbb P_\Bbb Cp^\perp$. If $\Gamma$ is a geodesic, then the point $p$ polar to the projective
line of $\Gamma$ belongs to $P$, a~contradiction. Therefore, there is a geodesic $G\subset P$ such that $G\cap\Gamma=\varnothing$,
i.e., $G\cap\Bbb P_\Bbb Cp^\perp=\varnothing$, implying that the bisector $B$ with the real spine $G$ does not contain $p$. As
$P\subset B$, we~have shown that $B_W=P$.

Conversely, let $P$ be an $\Bbb R$-plane. Then
$$W:=\big\{h\in{\Lin}_\Bbb C(V,V)\mid\tr h=0,\ h^*=h,\ \langle hp,p\rangle=0\text{ for all }p\in P\big\}$$
is an $\Bbb R$-linear subspace. In view of the above considerations, in order to show that $W$ is a linear family of bisectors, it
suffices to observe that $\det W=0$. Let $U\subset V$ be a totally real subspace such that $P=\Bbb P_\Bbb CU$ and let $h\in W$. From
$h^*=h$ and $\big\langle h(u_0+u_1),u_0+u_1\big\rangle=0$ for all $u_0,u_1\in U$, we infer $\Re\langle hu_0,u_1\rangle=0$. It follows
that $hU\subset iU$. Hence, there exists $\varphi\in\Lin_\Bbb R(U,U)$ such that $hu=i\varphi u$ for all $u\in U$. From
$\langle hu_0,u_1\rangle=\langle u_0,hu_1\rangle$ for $u_0,u_1\in U$, we conclude that
$\langle\varphi u_0,u_1\rangle+\langle u_0,\varphi u_1\rangle=0$. In~other words, $\varphi+\varphi^*=0$, where $\varphi^*$ stands for
the adjoint to $\varphi$ in the sense of the form $\langle-,-\rangle$ restricted on~$U$. It is well known (and can be easily verified
with a straightforward calculation) that all elements in the Lie algebra $\o U$ of the orthogonal group $\O U=\O(2,1)$ are degenerate.
Hence, $\varphi$ is degenerate and so is~$h$.

By Remark 1.2.1, no $3$ bisectors in $W$ are transversal at any $p\in P\setminus\S$
$_\blacksquare$

\medskip

{\bf4.4.2.~Lemma.} {\sl Let\/ $S$ be a projective line with the polar point\/ $p\notin\S$. Then
$$W_S:=\big\{h\in{\Lin}_\Bbb C(V,V)\mid\tr h=0,\ h^*=h,\ \langle hs,s\rangle=0\text{ \rm for all }s\in S\cup\{p\}\big\}$$
is a\/ $4$-dimensional\/ {\rm(}maximal\/{\rm)} nonconfocal family of bisectors and\/ $B_{W_S}=S\cup\{p\}$.

Conversely, any\/ $4$-dimensional\/ {\rm(}maximal\/{\rm)} nonconfocal family of bisectors has the above form.

No\/ $4$ bisectors of such a family are transversal at any $s\in S\setminus\S$. Any\/ $4$ bisectors that span the family are
transversal at\/ $p$.}

\medskip

{\bf Proof.} Let $U\subset V$ be a $2$-dimensional $\Bbb C$-linear subspace such that $S=\Bbb P_\Bbb CU$. Clearly, $W_S\subset W_U$.
By Corollary 3.5 and Remark 4.1, every $0\ne h\in W_U$ provides a bisector $B_h$ containing $S$. Therefore, $p$~belongs to the real
spine of $B_h$, hence, $p\in B_h$. Consequently, $p\in B_{W_U}$ and $W_S=W_U$. Since $S\cup\{p\}$ contains a negative point, $W_S$ is
a linear family of bisectors. It cannot be confocal because $\PU V$ acts transitively on points of a same signature in $S$. By
Proposition 1.2.3, we obtain an $\Bbb R$-linear embedding $f:W_S\hookrightarrow U$ such that $h(fh)=0$ for all $h\in W_S$.

Let $f_0\in S$ be nonisotropic. By Remark 2.5, we find a projective line $S'$ orthogonal to $S$ such that $f_0\in S'$. By Lemma 4.3.1,
(4.3.2) defines a confocal line $W$ of bisectors with the common focus $f_0$ and $S\subset B_W$. Hence, $W\subset W_S$ and
$fW\subset\Bbb Cf_0$. Therefore, $fW=\Bbb Cf_0$. Taking other nonisotropic points $f_0\in S$, we conclude that $fW_S=U$. In
particular, $\dim_\Bbb RW_S=4$ and $W_S$ is a maximal linear family of bisectors.

Let us show that $B_{W_S}=S\cup\{p\}$. Take any $q\notin S\cup\{p\}$ and pick $f_0\in S\setminus\S$ not belonging to the projective
line joining $p,q$. Denote by $L(f_0,q)$ the projective line joining $f_0,q$ and let $C:=\Bbb P_\Bbb Cf_0^\perp$. Then
$p\in C\not\ni f_0$ because $f_0\notin\S$. Since $f_0,p,q$ are not on a same projective line, we have $L(f_0,q)\cap C=\{d\}$ and
$d\ne p$. Pick a geodesic $R\subset C$ such that $p\in R\not\ni d$. Then $R$ is a real spine of some bisector $B_h$, $f_0$~is the
focus of $B_h$, and $S$ is a slice of $B_h$ because $p\in R$. So, $h\in W_U$. On the other hand, $q\notin B_h$ as, otherwise,
$L(f_0,q)\subset B_h$ and $d\in R$.

The converse follows from Proposition 1.2.3, Corollary 3.5, and the beginning of 4.4.

\smallskip

By Remark 1.2.1, no $4$ bisectors in $W_S$ are transversal at any $s\in S\setminus\S$ and any $4$ bisectors that span $W_S$ are
transversal at $p$
$_\blacksquare$

\medskip

{\bf4.4.3.~Lemma.} {\sl Let\/ $\Gamma$ be a geodesic and let\/ $p\in\Gamma\setminus\S$. Then
$$W:=\big\{h\in{\Lin}_\Bbb C(V,V)\mid\tr h=0,\ h^*=h,\ \langle hd,d\rangle=0\text{ \rm for all }d\in S\cup\Gamma\big\}.$$
is a\/ $3$-dimensional nonconfocal linear family of bisectors whose foci belong to the projective line\/ $S$ with the polar point\/
$p\notin\S$. We have\/ $B_W=S\cup\Gamma$.

Conversely, any\/ $3$-dimensional nonconfocal linear family of bisectors whose foci lie in a noneuclidean projective line has the
above form.

No\/ $3$ bisectors of such a family are transversal at any\/ $s\in S\setminus\S$. Any\/ $3$ bisectors that span the family are
transversal at any\/ $g\in\Gamma\setminus(\S\cup S)$.}

\medskip

{\bf Proof.} Let $W$ be a $3$-dimensional nonconfocal linear family of bisectors whose foci lie in a projective line $S$ with polar
point $p\not\in\S$ and let $U\subset V$ be a $2$-dimensional $\Bbb C$-linear subspace in $V$ such that $S=\Bbb P_\Bbb CU$. Then, by
Proposition 1.2.3 and the beginning of 4.4, $W\subset W_S$ and there is an $\Bbb R$-linear isomorphism $f:W_S\to U$ such that
$h(fh)=0$ for all $h\in W_S$.

Clearly, $fW\cap ifW=\Bbb Cf_0$ for some $f_0\in U$. So, $fW=\Bbb Cf_0+\Bbb Rq$. If $f_0\notin\S$, we choose $q$ orthogonal to $f_0$.
If $f_0\in\S$, then we can choose $q$ such that $f_0\ne q\in S\cap\S$ because $\Bbb P_\Bbb CfW=S$. Thus, we~assume that $f_0,p,q\in V$
is a $\Bbb C$-linear basis either orthogonal or with the Gram matrix $G:=\left[\smallmatrix0&0&1\\0&1&0\\1&0&0\endsmallmatrix\right]$.
Let $W':=f^{-1}(\Bbb Cf_0)$ and $h:=f^{-1}(q)$. By Lemma 4.3.1, $B_{W'}=S\cup S'$, where $S'$ is the projective line spanned by
$f_0,p$. Since $f_0\ne q\in S$, the focus $q$ of $B_h$ does not belong to $S'$ and, therefore, $B_h$ is the projective cone with the
apex $q$ and the base $\Gamma:=S'\cap B_h$. As $S$ is a slice of $B_h$ and $f_0\in S$, we have $f_0,p\in B_h$, hence,
$f_0,p\in\Gamma$. We claim that $\Gamma$ is a geodesic. Indeed, if $f_0,p,q$ are orthogonal, then $\Gamma$ is simply the real spine of
$B_h$. If $f_0,p,q$ have the Gram matrix $G$, then $B_h$ is parabolic, $q\ne f_0\in\S\cap B_h$, and $S'=\Bbb P_\Bbb Cf_0^\perp$. By
Remark 2.2, $\Gamma$ is an euclidean geodesic. It follows that $B_W=B_{W'}\cap B_h=(S\cup S')\cap B_h=S\cup\Gamma$.

Conversely, let $\Gamma$ be a geodesic and let $p\in\Gamma\setminus\S$. Then
$$W:=\big\{h\in{\Lin}_\Bbb C(V,V)\mid\tr h=0,\ h^*=h,\ \langle hd,d\rangle=0\text{ for all }d\in S\cup\Gamma\big\}\subset W_S,$$
where $S:=\Bbb P_\Bbb Cp^\perp$. We will show that $\dim_\Bbb RW=3$ and $B_W=S\cup\Gamma$.

Let $\{f_0\}:=S\cap\Gamma$, let $S'$ denote the projective line spanned by $f_0,p$ (so, $S$ and $S'$ are orthogonal), and let $W'$ be
given by (4.3.2). Then $\Gamma\subset S'$, $\dim_\Bbb RW'=2$, $B_{W'}=S\cup S'$, and $W'\subset W\subsetneq W_S$ by Lemma 4.3.1
because $B_W\supset\Gamma\not\subset B_{W_S}$.

Suppose that $f_0\notin\S$. Then we have an orthogonal $\Bbb C$-linear basis $f_0,p,q$, where $q\in S$ is the focus of the bisector
$B_h$ with the real spine $\Gamma$ and complex spine $S'$. It follows from $p\in\Gamma$ and $S'\not\subset B_h$ that $h\in W$ and
$h\notin W'$. Hence, $\dim_\Bbb RW_0=3$ and $B_W=B_{W'}\cap B_h=S\cup\Gamma$.

Suppose that $f_0\in\S$. Then $S$ is hyperbolic. Pick points $d,q$ such that $f_0,p\ne d\in\Gamma$ and $f_0\ne q\in\S\cap S$ and
consider the projective line $C:=\Bbb P_\Bbb Cq^\perp$ of signature $+0$. It follows from $q\notin\Gamma$ that $d\ne q$; denote by $b$
the point polar to the projective line $L(d,q)$ spanned by $d,q$. Then the points $b,p,q\in C$ are pairwise distinct because they are
the points polar to the projective lines $L(d,q),S,C$, each line contains $q$, which are distinct because their intersections with
$S'$, $q\notin S'$, are respectively $d,f_0,p\in\Gamma$. Since $S\cap C=\{q\}$ and $q\ne b\in C$, we see that $b\notin S$, implying
that the distinct points $b,p\in C$ are not orthogonal. Hence, there exists a unique geodesic $R$ such that $b,p\in R\subset C$. As
any geodesic in $C$ contains $q$, we obtain $b,p,q\in R$. Denote by $B_h$ the (parabolic) bisector with the real spine $R$. Then
$L(d,q),S,C$ are slices and $q$ is the focus of $B_h$. From $q\notin S'=\Bbb P_\Bbb Cf_0^\perp$ and $f_0\in S\subset B_h$, we conclude
that $f_0\in\S\cap B_h$ and, by Remark 2.2, that $S'\cap B_h$ is a geodesic. The geodesics $S'\cap B_h$ and $\Gamma$ coincide because
they have $3$ common points $d,f_0,p$. So, $h\in W\setminus W'$, $\dim_\Bbb RW_0=3$, and $B_W=B_{W'}\cap B_h=S\cup\Gamma$.

\smallskip

By Remark 1.2.1, no $3$ bisectors in $W_S$ are transversal at any $s\in S\setminus\S$ and any $3$ bisectors that span $W$ are
transversal at any $g\in\Gamma\setminus(\S\cup S)$
$_\blacksquare$

\medskip

{\bf4.4.4.~Lemma.} {\sl Let\/ $f_0\in \Gamma$ be a vertex of a hyperbolic geodesic\/ $\Gamma$. Then
$$W:=\big\{h\in{\Lin}_\Bbb C(V,V)\mid\tr h=0,\ h^*=h,\ \langle hd,d\rangle=0\text{ \rm for all }d\in\Gamma\cup S'\big\}$$
is a\/ $3$-dimensional nonconfocal linear family of bisectors whose foci belong to the projective line\/ $S'$ with the polar point\/
$f_0\in\S$. We have\/ $B_W=\Gamma\cup S'$.

Conversely, any\/ $3$-dimensional nonconfocal linear family of bisectors whose foci lie in an euclidean projective line has the above
form. Also, such families are exactly the\/ $1$-codimensional\/ $\Bbb R$-linear subspaces\/ $W\subset W_U$ subject to\/
$\langle-,f_0\rangle f_0\notin W$, where\/ $U:=f_0^\perp$ and\/ $f_0\in\S$.

No\/ $3$ bisectors of such a family are transversal at any\/ $s\in S'\setminus\{f_0\}$. Any\/ $3$ bisectors that span the family are
transversal at any\/ $g\in\Gamma\setminus\S$.}

\medskip

{\bf Proof.} Let $f_0\in\S$. Denote $U:=f_0^\perp$ and $S':=\Bbb P_\Bbb CU$. By Remark 4.1 and Corollary 3.5, $B_h$ is a bisector
containing $S'$ iff $h\in W_U\setminus\Bbb R\langle-,f_0\rangle f_0$. Given $p\in S'$, denote
$W_p:=\{h\in W_U\mid hp=0\}\ni\langle-,f_0\rangle f_0$. When $p\ne f_0$, a bisector $B_h$,
$h\in W_p\setminus\Bbb R\langle-,f_0\rangle f_0$, is determined by its real spine which can be any (hyperbolic) geodesic $R$ subject
to $f_0\in R\subset\Bbb P_\Bbb Cp^\perp$. For a fixed $p\ne f_0$, such $R$'s form a $1$-parameter family, hence, $\dim_\Bbb RW_p=2$.
When $p=f_0$, the real spine of $B_h$, $0\ne h\in W_{f_0}$, is an arbitrary geodesic $R\subset S'$. As the family of such geodesics is
$2$-dimensional, we conclude that $\dim_\Bbb RW_{f_0}=3$. Varying $p\in S'$, we arrive at $\dim_\Bbb RW_U=4$.

Let $W_0\subset W_U$ be a $3$-dimensional $\Bbb R$-linear subspace such that $W_0\cap\Bbb R\langle-,f_0\rangle f_0=0$. The bisectors
from $W_0$ cannot be all confocal because $W_0\not\subset W_p$ for any $p\in S'$ in view of dimension and
$\langle-,f_0\rangle f_0\notin W_0$. By Proposition 1.2.3, we obtain an $\Bbb R$-linear embedding $f:W_0\hookrightarrow V$ such that
$h(fh)=0$ for all $h\in W_0$. For some $p\in S'$, we have $fW_0\cap ifW_0\supset\Bbb Cp$, hence, $f^{-1}(\Bbb Cp)\subset W_p$. From
$f^{-1}(\Bbb Cp)\not\ni\langle-,f_0\rangle f_0\in W_p$ and $\dim_\Bbb Rf^{-1}(\Bbb Cp)=2$, we conclude that $p=f_0$. Therefore,
$\Bbb Cf_0\subset fW_0$ and $W':=f^{-1}(\Bbb Cf_0)$ is a confocal line of bisectors with the common focus $f_0$. By Lemma 4.3.1,
$B_{W'}=S\cup S'$, where $S$ is a projective line whose polar point $p$ belongs to $S'$, $f_0\ne p\in S'$, implying $f_0\in S\cap S'$.
From $\Bbb P_\Bbb CfW_0=S'$, we obtain $p=fh$ for some $h\in W_0\setminus W'$, hence, $W_0=\Bbb Rh+W'$. In other words, $p$ is the
focus of the (hyperbolic) bisector $B_h$, $S$ is the complex spine of $B_h$, $S'$ is a slice of $B_h$, and the hyperbolic geodesic
$\Gamma:=B_h\cap S\ni f_0$ is the real spine of $B_h$. Consequently, $B_{W_0}=B_h\cap B_{W'}=B_h\cap(S\cup S')=\Gamma\cup S'$ because
$S'\subset B_h$. It follows that $W_0$ is a maximal linear family of bisectors.

Let $W$ be a $3$-dimensional nonconfocal linear family $W$ of bisectors whose foci lie in a $S'$. Then, by~the beginning of 4.4,
$W\subset W_U$. By the above, $B_W=\Gamma\cup S'$, where $f_0\in\Gamma$ is a hyperbolic geodesic.

Conversely, let $f_0\in \Gamma$ be a vertex of a hyperbolic geodesic $\Gamma$. Put $U:=f_0^\perp$, $S':=\Bbb P_\Bbb CU$, and
$$W:=\big\{h\in{\Lin}_\Bbb C(V,V)\mid\tr h=0,\ h^*=h,\ \langle hd,d\rangle=0\text{ for all }d\in\Gamma\cup S'\big\}.$$
From $\Gamma\cap\B\ne\varnothing$, we infer that $W\subset W_U$ is a linear family of bisectors. By Remark 4.1 and Corollary 3.5,
$W\cap\Bbb R\langle-,f_0\rangle f_0=0$. Let $S$ denote the projective line of $\Gamma$. Then $S,S'$ are orthogonal and
$S\cap S'=\{f_0\}$. By Lemma 4.3.1, we obtain a confocal line $W'\subset W$ of bisectors with the common focus $f_0$. Thus,
$\dim_\Bbb RW=3$ because $h\in W\setminus W'$, where $B_h$ is the bisector with the real spine $\Gamma$.

\smallskip

By Remark 1.2.1, no $3$ bisectors in a linear family $W$ in question are transversal at any $s\in S'\setminus\{f_0\}$ and any $3$
bisectors that span $W$ are transversal at any $g\in\Gamma\setminus\S$
$_\blacksquare$

\medskip

{\bf4.5.~Nonconfocal lines of bisectors.} A nonconfocal $2$-dimensional linear family $W$ of bisectors is called a {\it nonconfocal
line\/} of bisectors. By Proposition 1.2.3, there exists a unique up to $\Bbb C^*$-proportionality $\Bbb R$-linear embedding
$f:W\hookrightarrow V$ such that $h(fh)=0$ for all $h\in W$ and $W\subset W_U$, where $U:=\Bbb CfW$ and $\dim_\Bbb CU=2$. By Corollary
3.5 and Remark 4.1, we get a common slice $S:=\Bbb P_\Bbb CU\subset B_W$ of the bisectors of $W$; denote by $p$ the point polar to
$S$. The circle $L:=\Bbb P_\Bbb CfW\subset S$ (independent of the choice of $f$) is called the {\it singular circle\/} of the
nonconfocal line of bisectors.

\medskip

{\bf4.5.1.~Definition.} Let $T\subset V$ be a $3$-dimensional $\Bbb R$-linear subspace such that $V=\Bbb CT$. The real projective
plane $Q:=\Bbb P_\Bbb CT$ is called a {\it geodesic cone\/} if $Q$ contains at least two geodesics. Let $p$ be the intersection of
geodesics $\Gamma_0,\Gamma_1\subset Q$. Then any real projective line $\Gamma'$ such that $p\in\Gamma'\subset Q$ is a geodesic.
Indeed, $\Gamma_i=\Bbb P_\Bbb CT_i$ and $T$ is the $\Bbb R$-span of $p,d_0,d_1\in V$, where $T_i$ stands for the $\Bbb R$-span of
$p,d_i$ and $\langle p,d_i\rangle\in\Bbb R$, $i=0,1$. Hence, $\Gamma'=\Bbb P_\Bbb CT'$, where $T'$ is the $\Bbb R$-span of
$p,d'$ and $d'$ belongs to the $\Bbb R$-span of $d_0,d_1$. Thus,~$\langle p,d'\rangle\in\Bbb R$ and $\Gamma'$ is a geodesic.

Of course, any $\Bbb R$-plane is a geodesic cone. If a geodesic cone $Q$ contains an extra geodesic $\Gamma$,
$p\notin\Gamma\subset Q$, then $Q$ is an $\Bbb R$-plane because we can take for $d_i$ the intersection $\Gamma\cap\Gamma_i$ and thus
obtain the $\Bbb R$-linear basis $p,d_0,d_1$ in $T$ with the real Gram matrix. We call $p$ the {\it apex\/} of the geodesic cone $Q$.
When $Q$ is an $\Bbb R$-plane, every $p\in Q$ serves as an apex of $Q$.

\medskip

{\bf4.5.2.~Lemma.} {\sl Let\/ $Q$ be a geodesic cone with apex\/ $p$ and let\/ $S$ be a projective line with the polar point\/ $p$.
Then
$$W:=\big\{h\in{\Lin}_\Bbb C(V,V)\mid\tr h=0,\ h^*=h,\ \langle hd,d\rangle=0\text{ \rm for all }d\in S\cap Q\big\}$$
is a nonconfical line of bisectors, $B_W=S\cup Q$, and the singular circle of $W$ equals $L:=S\cap Q$.

Conversely, any\/ nonconfocal line of bisectors has the above form.

No\/ $2$ bisectors of such a family are transversal at any\/ $l\in L\setminus\S$. Any\/ $2$ bisectors that span the family are
transversal at any\/ $d\in(S\cap Q)\setminus(L\cup\S)$.}

\medskip

{\bf Proof.} Let $W$ be a confocal line of bisectorts. By Lemmas 4.4.2, 4.4.3, and 4.4.4, there are distinct $3$-dimensional linear
families $W_0,W_1\subset W_U$ of bisectors such that $W=W_0\cap W_1$. It follows from $B_{W_i}=S\cup\Gamma_i$, $i=0,1$, that
$S\cup\Gamma_0\cup\Gamma_1\subset B_W$. The geodesics $\Gamma_i$'s are distinct because the $W_i$'s are distinct; both geodesics
contain the point polar to $S$.

It follows from $\dim_\Bbb RW=2$ and from the identity $h(fh)=0$, $h\in W$, that $W(fW)$ is a $1$-dimensional $\Bbb R$-linear space.
Hence, $W(fW)=\Bbb Rp$ for a suitable representative $p\in V$; this is the point polar to $S$ because $h(U)=\Bbb Cp$ for any
$0\ne h\in W$. Consequently, $p\in\Gamma_0\cap\Gamma_1$.

Pick $d_0\in\Gamma_0\setminus S$ such that $0\ne\langle p,d_0\rangle\in\Bbb Ri$ and let $T$ stand for the $\Bbb R$-span of $d_0,fW$.
Clearly, $Q:=\Bbb P_\Bbb CT$ is a real projective plane.

We claim that $B_W=S\cup Q$. Indeed, any point in $\Bbb P_\Bbb CV\setminus S$ has the form $d_0+fw+ifw'$ with $w,w'\in W$.
As $\langle hS',S'\rangle=0$ for every slice $S'$ of an arbitrary bisector $B_h$, we have $\langle hS,S\rangle=0$ for all $h\in W$.
Therefore, for any $h\in W$, we obtain $\big\langle h(fw+ifw'),fw+ifw'\big\rangle=0$ and, in view of $\langle hd_0,d_0\rangle=0$,
$$\big\langle h(d_0+fw+ifw'),d_0+fw+ifw'\big\rangle=\langle hd_0,fw+ifw'\rangle+\big\langle h(fw+ifw'),d_0\big\rangle=$$
$$=2\Re\big\langle h(fw+ifw'),d_0\big\rangle=2i\big\langle h(fw'),d_0\big\rangle$$
due to $h(fw)\in\Bbb Rp$ and $\langle p,d_0\rangle\in\Bbb Ri$. It follows that $d_0+fw+ifw'\in B_W$ iff $W(fw')=0$, i.e., iff $w'=0$.

Since $S\not\supset\Gamma_i\subset B_W$, we conclude that $\Gamma_i\subset Q$. In other words, $Q$ is a geodesic cone with the apex
$p$.

Clearly, $U\cap T\supset fW$ provides a real projective line inside $S\cap Q$. Since $S\cap Q=\Bbb P_\Bbb C(p^\perp\cap T)$ is a real
projective line, we see that $S\cap Q$ is the singular circle of $W$.

Conversely, let $Q$ be a geodesic cone with the apex $p$. Put $U:=q^\perp$, $S:=\Bbb P_\Bbb CU$, and
$$W:=\big\{h\in{\Lin}_\Bbb C(V,V)\mid\tr h=0,\ h^*=h,\ \langle hd,d\rangle=0\text{ for all }d\in S\cap Q\big\}.$$
Pick two distinct geodesics $p\in\Gamma_0,\Gamma_1\subset Q$ such that $\Gamma_0,\Gamma_1\not\subset S$. By Lemmas 4.4.2, 4.4.3, and
4.4.4, we~obtain distinct $3$-dimensional linear families $W_0,W_1\subset W_U$ of bisectors such that $B_{W_i}=S\cup\Gamma_i$,
$i=0,1$, and $W\subset W':=W_0\cap W_1$. Since $\Gamma_0\cup\Gamma_1\subset Q'$ and $B_{W'}=S\cup Q'$, where $Q'$ is a geodesic cone
with the apex $p$, we conclude that $Q'=Q$ and $W=W'$.

By Remark 1.2.1, no $2$ bisectors in $W$ are transversal at any $l\in L\setminus\S$ and any $2$ bisectors that span $W$ are
transversal at any $d\in(S\cup Q)\setminus(L\cup\S)$
$_\blacksquare$

\medskip

We warn the reader that, in general, the geodesics $\Gamma$ such that $p\in\Gamma\subset Q$ have nothing to do with the list of real
spines of the bisectors of the family $W$. Note also that it is quite possible that $L\cap(\B\cup\S)=\varnothing$; see, for example,
[AGG, p.~38, Criterion 4.3.3].

\smallskip

When $p\notin\S$, the $\Bbb R$-linear isomorphism $f:W_U\to U$ (unique up to $\Bbb C^*$-proportionality) provides the map
$\pi:\Gr_\Bbb R(2,W_U)\to\Gr_{\Bbb C|\Bbb R}(2,U)$ from the grassmannian $\Gr_\Bbb R(2,W_U)$ of $2$-dimensional $\Bbb R$-linear
subspaces in $W_U$ onto the $\Bbb C$-grassmannian $\Gr_{\Bbb C|\Bbb R}(2,U)$ of the circles in $S$ (see the end of [AGr, Examples~1.7]). The map $\pi$ sends $W\in\Gr_\Bbb R(2,W_U)$ onto the common focus $f_W\in S\subset\Gr_{\Bbb C|\Bbb R}(2,U)$ of $W$,
$\{f_W\}=\Bbb P_\Bbb C fW$, if $W$ is a confocal line, and onto the circle $L_W:=\Bbb P_\Bbb CfW\in\Gr_{\Bbb C|\Bbb R}(2,U)$,
otherwise. Outside the boundary $S$ of the $\Bbb C$-grassmannian $\Gr_{\Bbb C|\Bbb R}(2,U)$ the map $\pi$ is an $\U_1$-bundle because
we can describe the geodesic cone $Q$ such that $B_W=S\cup Q$ as $Q:=\Bbb P_\Bbb CT$, where $T:=\Bbb Rp+fW$ and $p\in V$ is a suitable
representative.

\medskip

{\bf4.5.3.~Corollary.} {\sl Any line of bisectors have a common slice}
$_\blacksquare$

\medskip

{\bf Proof of Proposition 1.2.4.} We just summarize Lemma 4.3.1, Remark 4.3.3, Lemmas 4.4.1, 4.4.2, 4.4.3, 4.4.4, and 4.5.2 and take
into account Remark 4.3.4
$_\blacksquare$

\medskip

{\bf4.6.~Lemma.} {\sl Let\/ $p\in B_{h_1}\cap B_{h_2}$ be a negative point different from the foci of bisectors\/ $B_{h_1},B_{h_2}$.
If\/~$B_{h_1},B_{h_2}$ are not transversal at\/ $p$, then\/ $B_{h_1},B_{h_2}$ belong to a nonconfocal line of bisectors and\/ $p$ lies
on the singular circle of the line, hence, on a common slice of\/ $B_{h_1},B_{h_2}$.}

\medskip

{\bf Proof.} We can assume that the $\Bbb R$-span $W$ of $h_1,h_2$ is $2$-dimensional. By Remark 1.2.1, $p$ is the focus of a
spherical bisector $B_h$, $0\ne h\in W$. We can assume that $h=h_2-h_1$. Hence, $h_1p=h_2p\in p^\perp$ and $h$ admits a
$\Bbb C$-linear basis of orthonormal eigenvectors $p,v_1,v_2$ with eigenvalues $0,\lambda,-\lambda$, $0\ne \lambda\in\Bbb R$. We can
choose representatives such that $h_1p=h_2p=a_1v_1+a_2v_2\ne0$ with $a_1,a_2\ge0$. In terms of the basis $p,v_1,v_2$, the condition
that $H^*=H$ for $H:=[h_{ij}]$ takes the form
$$h_{11}=\overline h_{11},\qquad h_{22}=\overline h_{22},\qquad h_{33}=\overline h_{33},\qquad h_{12}=-\overline h_{21},\qquad
h_{13}=-\overline h_{31},\qquad h_{23}=\overline h_{32}.$$
Taking into account that $\tr h_1=0$, we obtain $h_1=\left[\smallmatrix0&-a_1&-a_2\\a_1&r&\overline c\\a_2&c&-r\endsmallmatrix\right]$
and $h_2=\left[\smallmatrix0&-a_1&-a_2\\a_1&r+\lambda&\overline c\\a_2&c&-r-\lambda\endsmallmatrix\right]$, where $r\in\Bbb R$ and
$c\in\Bbb C$. Since $\det h_1=(a_2^2-a_1^2)r-2a_1a_2\Re c$, $\det h_2=(a_2^2-a_1^2)(r+\lambda)-2a_1a_2\Re c$, and $\lambda\ne0$,
we~conclude from $\det h_1=\det h_2=0$ that $a_1^2=a_2^2$ and $a_1a_2\Re c=0$. So, $a_1=a_2\ne0$ and $c\in\Bbb Ri$ because
$a_1v_1+a_2v_2\ne0$ and $a_1,a_2\ge0$. It follows that $\det W=0$. Therefore, $W$ is aline of bisectors.

By Remark 1.2.1, the only singular point in the base of a confocal line of bisectors is the common focus of the family. So, $W$ is not
confocal. The rest follows from Lemma 4.5.2
$_\blacksquare$

\medskip

{\bf4.7.~Lemma.} {\sl Let\/ $p\in B_{h_1}\cap B_{h_2}$ be a negative point different from the foci of\/ $B_{h_1},B_{h_2}$ that lies
neither in a common slice nor in a common meridian of\/ $B_{h_1},B_{h_2}$. If\/ $B_{h_1}\cap B_{h_2}\cap U\subset B_{h_3}$ for an open
neighbourhood\/ $U\ni p$, then\/ $h_3$ belongs to the\/ $\Bbb R$-span of\/ $h_1,h_2$.}

\medskip

{\bf Proof.} By Lemma 4.6, we can assume that $p\in D:=B_{h_1}\cap B_{h_2}\cap U$ is a smooth surface and that $U\subset\B$.

Suppose that the $\Bbb R$-span $W$ of $h_1,h_2,h_3$ is $3$-dimensional.

By Remark 1.2.1, every $d\in D$ is a focus of some bisector $B_{h_d}$, $h_d\in\Bbb P_\Bbb RW$. As $B_{h_1},B_{h_2}$ are transversal at
any $d\in D$, such a $h_d$ is unique by Remark 1.2.1. The algebraic map $\psi:D\to\Bbb P_\Bbb RW$ given by the rule
$\psi:d\mapsto h_d$ is injective because every $h_d$ has rank $2$. As the image of $\psi$ is Zariski dense, $\det W=0$. So, $W$ is a
linear family of bisectors. By Lemma 4.3.1, $B_{h_1},B_{h_2}$ cannot have a common focus since $p$ does not lie in a common slice of
$B_{h_1},B_{h_2}$. The family $W$ cannot be a one described in Lemma 4.4.3 or Lemma 4.4.4 because $p$ does not lie in a common slice
of $B_{h_1},B_{h_2}$ and $D\subset B_W$. Consequently, $W$ is an $\Bbb R$-plane of bisectors and $p\in D\subset B_W$ belongs to a
common meridian $B_W$ of $B_{h_1},B_{h_2}$ by Lemma 4.4.1 and Corollary 2.10, a contradiction
$_\blacksquare$

\medskip

{\bf4.8.~Corollary {\rm(G.~Giraud rigidity theorem)}.} {\sl Let\/ $B_{h_1},B_{h_2}$ be bisectors having no common slice and such
that\/ $D:=B_{h_1}\cap B_{h_2}\cap\B\ne\varnothing$. Then\/ $D$ is a smooth surface\/ {\rm(}possibly except at a single point, the
focus of some\/ $B_{h_i}${\rm)} and there exists at most one more bisector\/ $B_{h_3}$ that contains\/ a given nonempty open subset\/
$\varnothing\ne U\subset D$. Such a\/ $B_{h_3}$ contains\/ $B_{h_1}\cap B_{h_2}$.}

\medskip

{\bf Proof.} Since the bisectors $B_{h_1},B_{h_2}$ have no common slice, they cannot be confocal by Lemma 4.3.1 and $D$ cannot contain
both foci of $B_{h_1},B_{h_2}$. Suppose that $D=\{f_1\}$, where $f_1$ is the focus of $B_{h_1}$. Pick slices $S_i\subset B_{h_i}$,
$i=1,2$, such that $S_1\cap S_2=\{f_1\}$. As $f_1$ is not the focus of $B_{h_2}$, the intersection $S_1\cap S'_2$ is a negative point
different from $f_1$ for any slice $S'_2\subset B_{h_2}$ close to $S_2$, a contradiction.

By Lemma 4.6, $B_{h_1},B_{h_2}$ are transversal at any $p\in D\setminus\{f\}\ne\varnothing$, where $f$ is the focus of $B_{h_1}$ or of
$B_{h_2}$. So, $D\setminus\{f\}$ is a smooth surface. Denote by $W$ the $\Bbb R$-span of $h_1,h_2$. By Corollary 4.5.3, $W$~cannot be
a line of bisectors. In particular, by Lemma 4.4.1, $B_{h_1},B_{h_2}$ cannot have a common meridian. By~Lemma~4.7, $h_3\in W$. As $W$
is not a line of bisectors, the homogeneous polynomial $\det(h_1x_1+h_2x_2)$ of degree $3$ does not vanish identically
$_\blacksquare$

\bigskip

\centerline{\bf5.~Elliptic families of bisectors and equidistant loci}

\medskip

Besides the proofs of Lemma 1.2.6, Propositions 1.2.9 and 1.2.10, and Theorem 1.2.11, in this section, we find a criterion when a real
projective line study is tangent to an elliptic family $E_W$ of bisectors (Lemma 5.1), characterize singular points of a family $E_W$
equitant from $4$ points, and describe the zoology of equitant families in terms of the linear dependence between the points (Lemma
5.4).

\medskip

{\bf5.1.~Lemma.} {\sl Let\/ $B_h,B_{h'}$ be distinct bisectors such that\/ $\det W\ne0$, where\/ $W$ denotes the\/ $\Bbb R$-span of\/
$h,h'$. Then one of the following alternatives takes place.

\smallskip

$\bullet$ There exist pairwise distinct points\/ $p_1,p_2,p_3$ such that\/ $\langle-,p_i\rangle p_i-\langle-,p_{i+1}\rangle p_{i+1}$,
$i=1,2,3$ {\rm(}the~indices are modulo\/ $3${\rm),} are unique up to\/ $\Bbb R^*$-proportionality nonnull elements of rank\/ $\le 2$
in\/ $W$.

$\bullet$ The elements\/ $h,h'$ are unique up to\/ $\Bbb R^*$-proportionality nonnull elements of rank\/ $\le2$ in\/ $W$ and one of
the bisectors\/ $B_h,B_{h'}$ contains the focus of the other.

\smallskip

In the latter case, under the assumption that the focus of\/ $B_{h'}$ belongs to\/ $B_h$ {\rm(}which is equivalent to the assumption
that the real spine of\/ $B_h$ intersects the complex spine of\/ $B_{h'}${\rm)}, we have\/ $\det(hx+h'x')=rx^2x'$ with\/ $r\ne0$.}

\medskip

{\bf Proof.} The bisectors $B_h,B_{h'}$ cannot have a common focus because $\det W\ne0$. Hence, their complex spines $C,C'$ intersect
in a single point $\{p_2\}:=C\cap C'$. Denote by $R,R'$ the real spines of $B_h,B_{h'}$.

Suppose that $p_2\notin R\cup R'$. Reflecting $p_2$ in $R$ and in $R'$, we obtain by Lemma 2.12 the points $p_1\in C$ and $p_3\in C'$
such that $h=\langle-,p_1\rangle p_1-\langle-,p_2\rangle p_2$ and $h'=\langle-,p_2\rangle p_2-\langle-,p_3\rangle p_3$. Thus, we
arrive at the first alternative.

Since $p_2$ is the point polar to the projective line joining the foci of $B_h,B_{h'}$, the inclusion $p_2\in R\cup R'$ is equivalent
to the fact that one of the bisectors\/ $B_h,B_{h'}$ contains the focus of the other.

If $p_2\in R\cap R'$, then $B_h,B_{h'}$ have a common slice $S$ and, by Remark 4.1, $\det W=0$, a contradiction.

So, we can assume that $p_2\in R\setminus R'$, i.e., that $B_h$ contains the focus of $B_{h'}$. By Remark 2.3 and Lemma 2.12,
$h=\langle-,p_1\rangle ip_2-\langle-,p_2\rangle ip_1$ and $h'=\langle-,p_2\rangle p_2-\langle-,p_3\rangle p_3$, where $p_1\in R$ and
$p_1,p_2,p_3$ are not on a same projective line. In terms of the bases $q_1,q_2,q_3$ and $p_1,p_2,p_3$, where
$\big[\langle q_i,p_j\rangle\big]=\left[\smallmatrix1&0&0\\0&1&0\\0&0&1\endsmallmatrix\right]$, we have
$hx+h'x'=\left[\smallmatrix0&-ix&0\\ix&x'&0\\0&0&-x'\endsmallmatrix\right]$
$_\blacksquare$

\medskip

{\bf5.2.~Remark.} {\sl Let\/ $W$ be an elliptic family of bisectors. Any rank\/ $1$ point in\/ $\hat E_W$ is a complex point of a
confocal line included in\/ $E_W${\rm;} there are exactly two rank\/ $1$ points in\/ $\hat E_W$ if there is a confocal line inside\/
$E_W$.}

\medskip

{\bf Proof.} Since $h^*=h$ for any $h\in W$, it follows that $W\cap iW=0$. Let $h:=h_0+ih_1=\langle-,q_1\rangle q_0$ be a rank $1$
point in $\hat E_W$, $h_0,h_1\in W$. Then $2h_0=h^*+h=\langle-,q_0\rangle q_1+\langle-,q_1\rangle q_0$ and
$2h_1=i(h^*-h)=\langle-,q_0\rangle iq_1-\langle-,q_1\rangle iq_0$. Consequently, $h_0,h_1$ have rank $2$ and span a confocal line $U$
of bisectors included in $E_W$. As $\tr h=0$ implies $\langle q_0,q_1\rangle=0$, we obtain $B_U=S_0\cup S_1$ by Lemma 4.3.1, where
$S_0,S_1$ are orthogonal projective lines with the polar points $q_0,q_1$. There is at most one confocal line in $E_W$. Therefore,
$\langle-,q_0\rangle q_1$ and $\langle-,q_1\rangle q_0$ are unique rank $1$ points in $\hat E_W$
$_\blacksquare$

\medskip

{\bf Proof of Lemma 1.2.6.} Due to Remark 5.2, by sending $\hat h\in\hat E_W$ to its kernel, we obtain an algebraic morphism
$f:\hat E_W\to\Bbb P_\Bbb CV$. Let $\hat W:=W+iW\subset\Lin_\Bbb C(V,V)$ and let $\hat L_1\subset\Bbb P_\Bbb C\hat W$ be a comple
projective line such that $\hat L_1(e)=0$ for some $e\in\Bbb P_\Bbb CV$. Since $W$ contains no confocal line, $\hat L_1$ is not
defined over $\Bbb R$. Hence, $\hat E_W=\hat L\cup\hat L_1\cup\hat L_2$, where the complex projective line $\hat L_2$ is `conjugate'
to $\hat L_1$ and the complex projective line $\hat L$ is defined over $\Bbb R$. Denote by $L\subset E_W$ the real projective line
formed by the real points of $\hat L$. The point $\{p\}:=\hat L_1\cap\hat L_2$ is real and $\hat L_i$ cannot contain any other real
point. If $p\notin L$, then $E_W=L\cup\{p\}$ which contradicts the assumptions of Lemma~1.2.6. If~$p\in L$, then $E_W=L$ and there are
no $3$ noncollinear points in $E_W$.

Thus, the inverse algebraic morphism is well defined.
$_\blacksquare$

\medskip

{\bf5.3.~Lemma.} {\sl Let\/ $W$ be an equitant elliptic family of bisectors. Then\/ $W$ is real. Any real projective line\/
$L\subset E_W$ is spanned by a pair of bisectors given by\/ $w_1-w_2,w_3-w_4${\rm;} $w_1-w_3,w_2-w_4${\rm;} $w_1-w_4,w_2-w_3$. If\/
$E_W$ is reducible, then\/ $E_W$ is a smooth real conic plus a real projective line intersecting the conic in\/ $2$ points or\/ $E_W$
consists of\/ $3$ real projective lines sharing no common point.}

\medskip

{\bf Proof.} The listed pairs are the only pairs of the $6$ distinct bisectors given by $w_i-w_j$, $i\ne j$, that can belong to a real
projective line $L\subset E_W$ included in $E_W$ because, as observed in 1.2.8, $\sum_{i=1}^3r_iw_i$ (for example) has rank $3$ if
$\sum_{i=1}^3r_i=0$ and $r_i\ne0$ for all $i=1,2,3$. On the other hand, the $3$ distinct bisectors given by $w_1-w_2,w_2-w_3,w_3-w_1$
lie on a real projective line $L'$. Since $L\cap L'\ne\varnothing$, the bisector given, say, by $w_1-w_2$ belongs to $L$. For a
similar reason, one of the bisectors given by $w_2-w_3,w_3-w_4,w_4-w_2$ belongs to $L$. So, $L$ is spanned by the bisectors given by
$w_1-w_2,w_3-w_4$.

We conclude that $E_W$ can neither be a line plus a point, nor a double line plus a line. Let $L\subset E_W$ be a confocal line. No
pair can belong to $L$ because the $p_i$'s are not on a same complex projective line. Therefore, $W$ is generic.

The foci of the bisectors given by $w_1-w_2,w_1-w_3,w_1-w_4$ belong to the complex projective line $L$ with the
polar point $p_1$. If $L\subset E$, then $L=f\hat C$ for some irreducible component $\hat C$ of $\hat E_W$ containing the mentioned
$3$ bisectors. The foci of the other $3$ bisectors do not belong $L$ (if the focus $q$ of the bisector given, say, by $w_3-w_4$
belongs to $L$, then $q$ is orthogonal to $p_1,p_3,p_4$, a contradiction), and $\hat C$ cannot be a complex projective line. Hence,
the other $3$ bisectors have to belong to another component $\hat{C'}$ of $\hat E_W$ which has to be a complex projective line, a
contradiction. Consequently, $W$ is real.

$E_W$ cannot be $3$ real projective lines passing through a point because these $3$ lines are spanned by the mentioned $3$ pairs of
bisectors, hence, their common point $\sum_{i=1}^4r_iw_i$ satisfies $\sum_{i=1}^4r_i=0$ and $r_i+r_j=0$ for all $i\ne j$.

The remaining case when $E_W$ is a smooth real conic plus a real projective line is considered in the following lemma.

\medskip

{\bf5.4.~Lemma.} {\sl Let\/ $W$ be an equitant elliptic family of bisectors. A point\/ $h\in E_W$ is singular iff the focus of\/ $B_h$
belongs to\/ $B_W$.

Suppose that\/ $W$ is spanned by\/ $h_1,h_2,h_3\in E_W$ and let\/ $s\in B_W\setminus\S V$. Then the bisectors\/ $B_{h_i}$ are not
transversal at\/ $s$ iff\/ $s$ is the focus of a singular point of\/ $E_W$.

Keeping the\/ $w_i$'s, we can choose representative of the\/ $p_i$'s and their order so that\/ $a_1\ge a_2\ge a_3\ge a_4>0$ in a
essentially unique\/ $\Bbb C$-linear dependence\/ $\sum_{i=1}^4a_ip_i=0$ between the\/ $p_i$'s. Then one of the following
alternatives takes place.

\smallskip

$\bullet$ $a_1=a_2=a_3=a_4$, $E_W$ is formed by\/ $3$ real projective lines, the foci\/ $s_1,s_2,s_3$ of the pairwise intersections of
these lines are given by\/ $\langle s_1,p_1\rangle=\langle s_1,p_2\rangle=-\langle s_1,p_3\rangle=-\langle s_1,p_4\rangle$,
$\langle s_2,p_1\rangle=-\langle s_2,p_2\rangle=\langle s_2,p_3\rangle=-\langle s_2,p_4\rangle$,
$\langle s_3,p_1\rangle=-\langle s_3,p_2\rangle=-\langle s_3,p_3\rangle=\langle s_3,p_4\rangle$.

$\bullet$ $a_1=a_2>a_3=a_4$, $E_W$ is a smooth conic plus a real projective line transversal to the conic, the foci\/ $s_1,s_2$ of
the\/ $2$ intersection points of the conic and the line are given by\/
$\langle s_1,p_1\rangle=-\langle s_1,p_2\rangle=\langle s_1,p_3\rangle=-\langle s_1,p_4\rangle$ and\/
$\langle s_2,p_1\rangle=-\langle s_2,p_2\rangle=-\langle s_2,p_3\rangle=\langle s_2,p_4\rangle$.

$\bullet$ $a_1+a_4=a_2+a_3$, $a_1\ne a_2$, $E_W$ is an irreducible singular cubic, the focus\/ $s$ of the singular point of\/ $E_W$ is
given by\/ $\langle s,p_1\rangle=-\langle s,p_2\rangle=-\langle s,p_3\rangle=\langle s,p_4\rangle$.

$\bullet$ $a_1=a_2+a_3+a_4$, $E_W$ is an irreducible singular cubic, the focus\/ $s$ of the singular point of\/ $E_W$ is given by\/
$\langle s,p_1\rangle=-\langle s,p_2\rangle=-\langle s,p_3\rangle=-\langle s,p_4\rangle$.

$\bullet$ $a_1+a_4\ne a_2+a_3$, $a_1\ne a_2+a_3+a_4$, $E_W$ is a smooth cubic.}

\medskip

{\bf Proof.} The point $h$ is singular iff a generic real projective line passing through $h$ intersects $E_W$ in at most $2$ points
(not counting multiplicities). Thus, the first claim follows from Remark 5.1.

By Remark 1.2.1, the bisectors are not transversal at $s$ iff $s$ is a focus of some bisector $B_h$ of the family, implying the second
claim.

The fact that $s\in B_W$ is a focus of some bisector in $W$ means that
$\big|\langle s,p_1\rangle\big|=\big|\langle s,p_2\rangle\big|=\big|\langle s,p_3\rangle\big|=\big|\langle
s,p_4\rangle\big|$
and $\sum_{i=1}^4r_i\langle s,p_i\rangle p_i=0$ for suitable nontrivial $r_i$'s such that $\sum_{i=1}^4r_i=0$. Without loss of
generality, we can assume that $\langle s,p_1\rangle=1$ and $r_1=a_1$. Then the conditions take the form
$\langle s,p_i\rangle=:\varepsilon_i=\pm1$  with $r_i=\varepsilon_ia_i$ for all $i=2,3,4$ and $\sum_{i=1}^4\varepsilon_ia_i=0$. The
last equality guarantees the existence and uniqueness of $s\in V$ for given $\varepsilon_i$'s.

Let $L\subset E_W$ be a real projective line containing the bisectors $B_1$ and $B_2$ given by $w_i-w_j$ and $w_k-w_l$. As $B_1,B_2$
have a common slice spanned by their foci $q_1\ne q_2$, we obtain
$\big|\langle q_2,p_i\rangle\big|=\big|\langle q_2,p_j\rangle\big|\ne0$. From $a_i\langle q_2,p_i\rangle+a_j\langle q_2,p_j\rangle=0$,
we conclude that $a_i=a_j$ and, by symmetry, $a_k=a_l$.

Now we can see that the cases when $E_W$ contains a real projective line are listed in the first $2$ alternatives.

Taking into account that the singular point of an irreducible singular real cubic is real, we arrive at the remaining $3$
alternatives
$_\blacksquare$
$_\blacksquare$

\medskip

{\bf Proof of Proposition 1.2.9.} Let $W$ be a real elliptic family of bisectors and let $L\subset\Bbb P_\Bbb CV$, $L\not\subset E$,
be a complex projective line containing the foci of $3$ distinct bisectors $h_1,h_2,h_3\in E_W$. The point $p$ polar to $L$ is a
common point of the complex spines $C_1,C_2,C_3$ of $B_{h_1},B_{h_2},B_{h_3}$. Denote by $R_i\subset C_i$ the real spine of $B_{h_i}$.
Since $f:\hat E_W\to E\subset\Bbb P_\Bbb CV$ is an isomorphism, the $C_i$'s are pairwise distinct.

Suppose that $p\in R_1\cap R_2\cap R_3$. Then the bisectors $B_{h_1},B_{h_2},B_{h_3}$ have a common slice. If $h_1,h_2,h_3$ are not
collinear, it follows from (3.6) and Remark 4.1 that $\det W=0$, a contradiction. If all $3$ belong to a nonconfocal real projective
line $L'\subset E_W$, the component $f\hat{L'}$ of $E$ has at least $3$ common points with $L$. As $f\hat{L'}\ne E$, we conclude that
$E\supset f\hat{L'}=L$, a contradiction.

Suppose that $R_1\cap R_2\ni p\notin R_3$. Then the bisectors $B_{h_1},B_{h_2}$ generate a nonconfocal real projective line
$L'\subset E_W$ and $fL'\subset L$ by the beginning of 4.5, hence, $f\hat{L'}=L$, implying $L\subset E$, a contradiction.

Suppose that $p\notin R_1\cup R_2\cup R_3$. By reflecting $p_4:=p$ in $R_i$ (see Lemma 2.12), we get $p_i\in C_i$ such that
$h_i=\langle-,p_i\rangle p_i-\langle-,p_4\rangle p_4$ for all $i=1,2,3$. The only $3$ points among the $p_i$'s that can belong to a
same complex projective line are $p_2,p_3,p_4$. In this case, the bisectors $B_{h_1-h_2},B_{h_2-h_3},B_{h_3-h_1}$ are confocal. So,
they have to coincide. In the pairs of points $p_2,p_3$, $p_3,p_4$, $p_4,p_2$, one point is the other reflected in the real spine
of this bisector. This implies $p_2=p_3=p_4$, a contradiction. Consequently, $W$~is equitant from the $p_i$'s.

We will reduce the remaining case $R_3\ni p\notin R_1\cup R_2$, to the previous one. Since the points $fh_1,fh_2,fh_3\in E\cap L$ are
pairwise distinct, they are smooth in $E$, hence, the points $h_1,h_2,h_3\in E_W$ are smooth in $E_W$. Denote by $L_i$ the real
projective line joining $h_i$ and $h_3$, $i=1,2$, and by $L_3$, the real projective line joining $h_1$ and $h_2$. By Lemma 5.1, the
relation $R_3\ni p\notin R_1\cup R_2$ means exactly that $L_i\cap E_W=2h_i+h_3$ for all $i=1,2$, i.e., that $L_i$ is tangent to $E_W$
at $h_i$ and $L_i\not\subset E_W$

Suppose that $E_W$ is irreducible. Then, using the standard arguments involving linear systems of divisors, we will `move'
$h_1,h_2,h_3$ so that the relation of the type $R_3\ni p\notin R_1\cup R_2$ will disappear. As $E_W$ is infinite and, for a given
smooth point $h\in E_W$, there are at most $4$ nonsingular points $h'\in E_W\setminus\{h\}$ such that the real projective line joining
$h$ and $h'$ is tangent to $E_W$ at $h'$, we can pick a smooth point $h'_3\in E_W\setminus\big(L_3\cup\{h_3\}\big)$ such that the real
projective line $L'_i$ joining $h_i\in L_3$ and $h'_3$ is not tangent to $E_W$ at $h'_3$ for all $i=1,2$. As $h'_3\ne h_3$, the line
$L'_i$ is not tangent to $E_W$ at any point. Hence, $L'_2\cap E_W=h_2+h'_2+h'_3$ for a suitable smooth point $h'_2\in E_W$, where
$h_2,h'_2,h'_3$ are distinct. The points $h_1,h'_2,h'_3$ could not be on a same real projective line because such a line should
coincide with $L'_2$ and then with $L_3$, which would contradict $h'_3\notin L_3$. Since $h_1,h'_3\in L'_1$, $h'_2,h'_3\in L'_2$, and
$L'_1,L'_2$ are not tangent to $E_W$ at any point, it remains to show that $fh_1$ lies on the complex projective line
$L'\subset\Bbb P_\Bbb CV$ joining $fh'_2,fh'_3$.

Denote by $l_2,l'_2$ the linear forms providing the lines $\hat L_2,\hat L'_2$ and by $l,l'$, the linear forms providing the lines
$L,L'\subset\Bbb P_\Bbb CV$. Clearly, $\varphi:=l_2/l'_2$ and $\psi:=l/l'$ are nonnull rational functions respectively on $\hat E_W$
and $E$ that are regular and invertible at the singular point, if it exists. For a suitable smooth point $q\in L'\cap E$, the divisors
of $\varphi$ and $\psi$ are $h_2+h_3-h'_2-h'_3$ and $fh_1+fh_2+fh_3-q-fh'_2-fh'_3$, hence, $q-fh_1$ is the divisor of
$f_*\varphi/\psi$. For an irreducible cubic, this is well known to imply $q=fh_1$, and we are done.

Suppose that $E_W$ is reducible. Since $E_W\not\supset L_i$ is tangent to $E_W$ for $i=1,2$, the cubic $E_W$ should consist of a real
smooth conic $C$ and a real projective line $R$, $E_W=C\cup R$, so that $h_1,h_2\in C\setminus R$ and $h_3\in R\setminus C$. As above,
$h_i=\langle-,p_i\rangle p_i-\langle-,p_3\rangle p_3$ for $i=1,2$ by Lemma 2.12, where $p_3:=p$, and
$h_3=\langle-,p_3\rangle ip_4-\langle-,p_4\rangle ip_3$ by Remark 2.3, where $p_3,p_4\in R_3$ and $\langle p_3,p_4\rangle\in\Bbb R$.
Clearly, $p_1,p_2,p_3$ are $\Bbb C$-linearly independent. We can choose representatives of $p_1,p_2$ such that
$p_4=\sum_{i=1}^3a_ip_i$ with $a_1,a_2\ge0$; denote $a:=\Im a_3$. We define the basis $q_1,q_2,q_3$ in $V$ by means of
$\big[\langle q_i,p_j\rangle\big]=\left[\smallmatrix1&0&0\\0&1&0\\0&0&1\endsmallmatrix\right]$. In~terms of the bases $q_1,q_2,q_3$
and $p_1,p_2,p_3$, we have
$x_1h_1x+x_2h_2+x_3h_3=\left[\smallmatrix x_1&0&ia_1x_3\\0&x_2&ia_2x_3\\-ia_1x_3&-ia_2x_3&-x_1-x_2-2ax_3\endsmallmatrix\right]$.
Consequently, $E_W$ is given by the equation $x_1x_2(x_1+x_2+2ax_3)+(a_1^2x_2+a_2^2x_1)x_3^2=0$ and $R$, by the equation of the form
$b_1x_2+b_2x_1=0$. If $a_1=a_2=0$, $E_W$ would be formed by $3$ real projective lines. Hence, we can assume that $b_i:=a_i^2$ and
$x_1x_2(x_1+x_2+2ax_3)$ should be divisible by $a_1^2x_2+a_2^2x_1$. Up to symmetry, either $a_1=a_2\ne0=a$ or $a_1\ne0=a_2$. In the
second case, the point $h_1$ does not satisfy the equation $x_1(x_1+x_2+2ax_3)+a_1^2x_3^2=0$ of $C$. In the first case, we can take a
suitable representative of $p_4$ providing $a_1=a_2=1$. Thus, $C$ is given by $x_1x_2+x_3^2=0$ and $R$, by $x_1+x_2=0$.
The points $q_2,2q_1-2q_2+iq_3,4q_1-q_2+2iq_3$ are the foci of the bisectors given by $h_1$, $h'_2:=h_1-h_2+2h_3$,
$h'_3:=h_1-4h_2+2h_3$. These foci lie on a same complex projective line not included in $E$. The points $h_2,h'_2,h'_3\in E_W$ are
collinear as well as the points $h_1,2h_1-2h_2+h_3,h'_3\in E_W$.

The converse follows from Lemma 5.3
$_\blacksquare$

\medskip

Let $W$ be an elliptic family of bisectors equitant from normalized points $p_i\in V$, $1\le i\le4$, of signature $\sigma$,
$\langle p_i,p_i\rangle=\sigma$. As in Lemma 5.4, we assume that $a_1\ge a_2\ge a_3\ge a_4>0$ in an essentially unique
$\Bbb C$-linear dependence $\sum_{i=1}^4a_ip_i=0$ between the $p_i$'s.

\medskip

{\bf Proof of Proposition 1.2.10.} Since
$B_W=\Big\{q\mid\big|\langle q,p_1\rangle\big|=\big|\langle q,p_2\rangle\big|=\big|\langle q,p_3\rangle\big|=\big|\langle
q,p_4\rangle\big|\Big\}$,
in terms of the coordinates $x_i:=\langle-,p_i\rangle$, $i=1,2,3$, we have
$$B_W=\big\{[x_1,x_2,1]\in\Bbb P_\Bbb CV\mid|x_1|=|x_2|=1,\ |a_1x_1+a_2x_2+a_3|=a_4\big\}.$$
Let $x_1:=\frac{(s_0+is_1)^2}{s_0^2+s_1^2}$ and $x_2:=\frac{(t_0+it_1)^2}{t_0^2+t_1^2}$ with $s_0,s_1,t_0,t_1\in\Bbb R$, thus
providing $|x_1|=|x_2|=1$. Then the condition $|a_1x_1+a_2x_2+a_3|^2=a_4^2$ means that
$$(a_1^2+a_2^2+a_3^2-a_4^2)(s_0^2+s_1^2)(t_0^2+t_1^2)+2a_1a_2\big((s_0t_0+s_1t_1)^2-(s_0t_1-s_1t_0)^2\big)+$$
$$+2a_1a_3(s_0^2-s_1^2)(t_0^2+t_1^2)+2a_2a_3(t_0^2-t_1^2)(s_0^2+s_1^2)=0,$$
i.e., $p(s_0,s_1,t_0,t_1)=0$, where
$$p(s_0,s_1,t_0,t_1):=r_{00}s_0^2t_0^2+r_{01}s_0^2t_1^2+2s_0s_1t_0t_1+r_{10}s_1^2t_0^2+r_{11}s_1^2t_1^2$$
and
$$r_{00}:=\frac{(a_1+a_2+a_3+a_4)(a_1+a_2+a_3-a_4)}{4a_1a_2},\quad r_{01}:=\frac{(a_1-a_2+a_3+a_4)(a_1-a_2+a_3-a_4)}{4a_1a_2},$$
$$r_{10}:=\frac{(a_1-a_2-a_3+a_4)(a_1-a_2-a_3-a_4)}{4a_1a_2},\quad r_{11}:=\frac{(a_1+a_2-a_3+a_4)(a_1+a_2-a_3-a_4)}{4a_1a_2}.$$

First, we are going to understand when $p(s_0,s_1,t_0,t_1)$ has a linear factor. If $p(s_0,s_1,t_0,t_1)$ is divisible by
$c_0s_0+c_1s_1$, then
$$p(s_0,s_1,t_0,t_1)=(c_0s_0+c_1s_1)\big(s_0q_0(t_0,t_1)+s_1q_1(t_0,t_1)\big)=$$
$$=s_0s_1\big(c_0q_1(t_0,t_1)+c_1q_0(t_0,t_1)\big)+c_0s_0^2q_0(t_0,t_1)+c_1s_1^2q_1(t_0,t_1),$$
implying
$$c_0q_0(t_0,t_1)=r_{00}t_0^2+r_{01}t_1^2,\qquad c_1q_1(t_0,t_1)=r_{10}t_0^2+r_{11}t_1^2,\qquad
c_0q_1(t_0,t_1)+c_1q_0(t_0,t_1)=2t_0t_1.$$
If $q_0(t_0,t_1)$ has the term $t_0t_1$, then $c_0=0$ and $p(s_0,s_1,t_0,t_1)$ is divisible by $s_1$. Otherwise, $q_1(t_0,t_1)$ has
the term $t_0t_1$, implying that $c_1=0$ and that $p(s_0,s_1,t_0,t_1)$ is divisible by $s_0$. By symmetry between $s_0,s_1$ and
$t_0,t_1$, we conclude that only $s_0,s_1,t_0,t_1$ can be linear divisors of $p(s_0,s_1,t_0,t_1)$, providing the following cases:
$$r_{00}=r_{01}=0,\qquad r_{00}=r_{10}=0,\qquad r_{01}=r_{11}=0,\qquad r_{10}=r_{11}=0.$$
The equality $r_{00}=0$ is impossible in view of $a_1\ge a_2\ge a_3\ge a_4>0$. If $r_{11}=0$, then $a_1=a_2=a_3=a_4$ and, by Lemma
5.4, the family is formed by $3$ lines, i.e., $E_W$ is reducible.

Assume now that $p(s_0,s_1,t_0,t_1)$ has no linear factor. A decomposition of $p(s_0,s_1,t_0,t_1)$ of the bihomogeneous type
$(2,0)+(0,2)$ would imply a factor of degree $1$. So, if $p(s_0,s_1,t_0,t_1)$ is reducible, then
$p(s_0,s_1,t_0,t_1)=\big(s_0p_0(t_0,t_1)+s_1p_1(t_0,t_1)\big)\big(s_0q_0(t_0,t_1)+s_1q_1(t_0,t_1)\big)$ and
$$p_0(t_0,t_1)q_0(t_0,t_1)=r_{00}t_0^2+r_{01}t_1^2,\qquad p_1(t_0,t_1)q_1(t_0,t_1)=r_{10}t_0^2+r_{11}t_1^2,$$
$$p_0(t_0,t_1)q_1(t_0,t_1)+p_1(t_0,t_1)q_0(t_0,t_1)=2t_0t_1.$$
Since $p_0(t_0,t_1),p_1(t_0,t_1),q_0(t_0,t_1),q_1(t_0,t_1)\ne0$, we obtain
$$p_i(t_0,t_1)=\alpha_i(\beta_{i0}t_0+\beta_{i1}t_1),\qquad q_i(t_0,t_1)=\alpha_i^{-1}(\beta_{i0}t_0-\beta_{i1}t_1)$$
for suitable $\alpha_i\in\Bbb C$ and $\beta_{ij}\in\Bbb R\cup i\Bbb R$ such that $\beta_{i0}^2=r_{i0}$ and $\beta_{i1}^2=-r_{i1}$.
Thus,
$$(\alpha_0\alpha_1^{-1}+\alpha_0^{-1}\alpha_1)\beta_{00}\beta_{10}=0,\qquad(\alpha_0\alpha_1^{-1}+
\alpha_0^{-1}\alpha_1)\beta_{01}\beta_{11}=0,$$
$$(\alpha_0^{-1}\alpha_1-\alpha_0\alpha_1^{-1})(\beta_{00}\beta_{11}-\beta_{01}\beta_{10})=2.$$

If $\alpha_0\alpha_1^{-1}+\alpha_0^{-1}\alpha_1\ne0$, then $r_{01}=0$, i.e., $a_1=a_2$ and $a_3=a_4$. By Lemma 5.4, $E_W$ is
reducible. Therefore, $\alpha_0\alpha_1^{-1}+\alpha_0^{-1}\alpha_1=0$. It follows that $\alpha_1=\pm i\alpha_0$ and
$\beta_{00}\beta_{11}-\beta_{01}\beta_{10}=\mp i$. Without loss of generality, we take $\alpha_0=1$ and $\alpha_1=i$. Then
$$p(s_0,s_1,t_0,t_1)=d(s_0,s_1,t_0,t_1)q(s_0,s_1,t_0,t_1),$$
where
$$d(s_0,s_1,t_0,t_1):=\beta_{00}s_0t_0+\beta_{01}s_0t_1+i\beta_{10}s_1t_0+i\beta_{11}s_1t_1,$$
$$q(s_0,s_1,t_0,t_1):=\beta_{00}s_0t_0-\beta_{01}s_0t_1-i\beta_{10}s_1t_0+i\beta_{11}s_1t_1.$$
Since $q(s_0,s_1,t_0,t_1)=d(s_0,-s_1,t_0,-t_1)$, we study only zeros of $d(s_0,s_1,t_0,t_1)$.

The equation $d(s_0,s_1,t_0,t_1)=0$ can be written as
$\left[\smallmatrix s_0&s_1\endsmallmatrix\right]\left[\smallmatrix
a_{00}&a_{01}\\a_{10}&a_{11}\endsmallmatrix\right]\left[\smallmatrix t_0\\t_1\endsmallmatrix\right]=0$,
where $a_{0i}:=\beta_{0i}$ and $a_{1i}:=i\beta_{1i}$. It follows from $\beta_{00}\beta_{11}-\beta_{01}\beta_{10}=-i$ that
$\det\left[\smallmatrix a_{00}&a_{01}\\a_{10}&a_{11}\endsmallmatrix\right]=1$.

Suppose that the equation $d(s_0,s_1,t_0,t_1)=0$ admits infinitely many solutions in $\Bbb P_\Bbb R^1\times\Bbb P_\Bbb R^1$. Since,
for a given $\left[\smallmatrix s_0&s_1\endsmallmatrix\right]$, there is an essentially unique
$\left[\smallmatrix t_0&t_1\endsmallmatrix\right]$ subject to the equation, we only need to study if
$\left[\smallmatrix s_0&s_1\endsmallmatrix\right]\left[\smallmatrix a_{00}&a_{01}\\a_{10}&a_{11}\endsmallmatrix\right]$ is
proportional to a real one for infinitely many (hence, for all) nonproportional real
$\left[\smallmatrix s_0&s_1\endsmallmatrix\right]$. (In particular, this implies that
$\left[\smallmatrix a_{00}&a_{01}\endsmallmatrix\right]$ and $\left[\smallmatrix a_{10}&a_{11}\endsmallmatrix\right]$ are proportional
to real ones.) Taking $a_{ij}\in\Bbb R\cup i\Bbb R$ into account, we conclude that the $a_{ij}$'s are all real or they are all purely
imaginary. In~other words, $\pm p(s_0,s_1,t_0,t_1)=(b_{00}s_0t_0+b_{11}s_1t_1)^2-(b_{01}s_0t_1+b_{10}s_1t_0)^2$ with
$b_{ij}\in\Bbb R$. Since $r_{00}>0$, the sign in $\pm p(s_0,s_1,t_0,t_1)$ is $+$. So, $r_{01}=-b_{01}^2\le0$, a contradiction.

By Lemma 5.4, any point $p\in B_W\cap\B V$ different from the focus of a singular point of $E_W$ provides infinitely many points in
$B_W$ close to $p$. So, if $p(s_0,s_1,t_0,t_1)=0$ admits only finitely many solutions in $\Bbb P_\Bbb R^1\times\Bbb P_\Bbb R^1$, then
$a_1+a_4=a_2+a_3$ or $a_1=a_2+a_3+a_4$ by Lemma 5.4. Hence, $\beta_{10}=0$, and it follows from
$\beta_{00}\beta_{11}-\beta_{01}\beta_{10}=-i$ that $r_{00}r_{11}=1$. If $a_4=a_2+a_3-a_1$, then
$r_{00}r_{11}=\frac{(a_2+a_3)(a_1-a_3)}{a_1a_2}$. Therefore, $(a_1-a_2-a_3)a_3=0$, implying $a_1-a_2-a_3=0$ and $a_4=0$, a
contradiction. If $a_4=a_1-a_2-a_3$, then $r_{00}r_{11}=\frac{(a_2+a_3)(a_1-a_3)}{a_1a_2}$. Consequently, $(a_1-a_2-a_3)a_3=0$ and
again $a_4=0$, a contradiction.
$_\blacksquare$

\medskip

{\bf Proof of Theorem 1.2.11.} Take $h\in\Lin_\Bbb C(V,V)$ such that $h^*=h$ and $\langle hb,b\rangle=0$ for all $b\in B$. In terms of
the coordinates $x_i$'s introduced at the beginning of the proof of Proposition 1.2.10, the hermitian form
$[u,v]:=\langle hu,v\rangle$ has a hermitian matrix $[h_{jk}]$. By Proposition 1.2.10, $x_1:=\frac{(s_0+is_1)^2}{s_0^2+s_1^2}$ and
$x_2:=\frac{(t_0+it_1)^2}{t_0^2+t_1^2}$ subject to $p(s_0,s_1,t_0,t_1)=0$ provide an irreducible curve. Since $B$ is Zariski dense in
this curve, the~equality $p(s_0,s_1,t_0,t_1)=0$ implies
$h_{11}+h_{22}+h_{33}+2\Re(h_{12}x_1\overline x_2+h_{23}x_2+h_{31}\overline x_1)=0$. In terms of $s_0,s_1,t_0,t_1$, it takes the form
$$(h_{11}+h_{22}+h_{33})(s_0^2+s_1^2)(t_0^2+t_1^2)+2h_{12}^0\big((s_0t_0+s_1t_1)^2-(s_0t_1-s_1t_0)^2\big)+
4h_{12}^1(s_0t_0+s_1t_1)(s_0t_1-s_1t_0)+$$
$$+2h_{23}^0(s_0^2+s_1^2)(t_0^2-t_1^2)-4h_{23}^1(s_0^2+s_1^2)t_0t_1+2h_{31}^0(s_0^2-s_1^2)(t_0^2+t_1^2)+
4h_{31}^1s_0s_1(t_0^2+t_1^2)=0,$$
where $h_{jk}=h_{jk}^0+ih_{jk}^1$ with $h_{jk}^0,h_{jk}^1\in\Bbb R$. Hence,
$$(h_{11}+h_{22}+h_{33}+2h_{12}^0+2h_{23}^0+2h_{31}^0)s_0^2t_0^2+(h_{11}+h_{22}+h_{33}-2h_{12}^0-2h_{23}^0+2h_{31}^0)s_0^2t_1^2+
8h_{12}^0s_0s_1t_0t_1+$$
$$+(h_{11}+h_{22}+h_{33}-2h_{12}^0+2h_{23}^0-2h_{31}^0)s_1^2t_0^2+(h_{11}+h_{22}+h_{33}+2h_{12}^0-2h_{23}^0-2h_{31}^0)s_1^2t_1^2+$$
$$+4(h_{12}^1-h_{23}^1)s_0^2t_0t_1+4(h_{31}^1-h_{12}^1)s_0s_1t_0^2+4(h_{12}^1+h_{31}^1)s_0s_1t_1^2-
4(h_{12}^1+h_{23}^1)s_1^2t_0t_1=0.$$
By Proposition 1.2.10, the left-hand side should be proportional to $p(s_0,s_1,t_0,t_1)$. In particular,
$h_{12}^1=h_{23}^1=h_{31}^1=0$.

If $h_{12}^0=0$, then
$$h_{11}+h_{22}+h_{33}+2h_{23}^0+2h_{31}^0=0,\qquad h_{11}+h_{22}+h_{33}-2h_{23}^0+2h_{31}^0=0,$$
$$h_{11}+h_{22}+h_{33}+2h_{23}^0-2h_{31}^0=0,\qquad h_{11}+h_{22}+h_{33}-2h_{23}^0-2h_{31}^0=0,$$
implying $h_{23}^0=h_{31}^0=0$ and $h_{11}+h_{22}+h_{33}=0$. In this case $h$ belongs to the $\Bbb R$-span of $w_1-w_2$ and $w_2-w_3$.

If $h_{12}^0\ne 0$, we can assume that the left-hand side is equal to $4a_1a_2p(s_0,s_1,t_0,t_1)$. Hence,
$$h_{11}+h_{22}+h_{33}+2h_{12}^0+2h_{23}^0+2h_{31}^0=(a_1+a_2+a_3+a_4)(a_1+a_2+a_3-a_4),$$
$$h_{11}+h_{22}+h_{33}-2h_{12}^0-2h_{23}^0+2h_{31}^0=(a_1-a_2+a_3+a_4)(a_1-a_2+a_3-a_4),$$
$$h_{11}+h_{22}+h_{33}-2h_{12}^0+2h_{23}^0-2h_{31}^0=(a_1-a_2-a_3+a_4)(a_1-a_2-a_3-a_4),$$
$$h_{11}+h_{22}+h_{33}+2h_{12}^0-2h_{23}^0-2h_{31}^0=(a_1+a_2-a_3+a_4)(a_1+a_2-a_3-a_4),$$
$$h_{12}^1=h_{23}^1=h_{31}^1=0,\qquad h_{12}^0=a_1a_2.$$
Subtracting the second equality from the first one and taking into account the last one, we obtain $h_{23}^0=a_2a_3$. Subtracting the
second equality from the fourth one and taking into account the last one, we obtain $h_{31}^0=a_3a_1$. Now it is clear that
$h_{11}+h_{22}+h_{33}=a_1^2+a_2^2+a_3^2-a_4^2$ and that
$[h_{jk}]=\left[\smallmatrix h_{11}&a_1a_2&a_1a_3\\a_2a_1&h_{22}&a_2a_3\\a_3a_1&a_3a_2&h_{33}\endsmallmatrix\right]$.
Without loss of generality, we can assume that $a_4=1$. Then the matrix of the hermitian form
$[u,v]_4:=\langle w_4u,v\rangle=\langle u,p_4\rangle\langle p_4,v\rangle$ equals
$\left[\smallmatrix a_1^2&a_1a_2&a_1a_3\\a_2a_1&a_2^2&a_2a_3\\a_3a_1&a_3a_2&a_3^2\endsmallmatrix\right]$. So,
$h=w_4+\sum_{i=1}^3(h_{ii}-a_i^2)w_i$. It remains to observe that $1+\sum_{i=1}^3(h_{ii}-a_i^2)=0$
$_\blacksquare$

\bigskip

\centerline{\bf6.~Appendix: Proof of Theorem 1.3.2}

\medskip

{\bf6.1.~Lemma.} {\sl Let\/ $L,M\subset\Bbb P_\Bbb KA$ be lines, let\/ $p_1,p_2,p_3\in L$, and let\/ $q_1,q_2,q_3\in M$. Suppose that
no other\/ $3$ points among\/ $p_1,p_2,p_3,q_1,q_2,q_3$ lie on a line and that neither\/ $p'_1,p'_2,p'_3\in\Bbb P_\Bbb KA$ nor\/
$q'_1,q'_2,q'_3\in\Bbb P_\Bbb KA$ lie on a line. Then the points\/ $p_i\otimes p'_i$, $q_j\otimes q'_j$, $i,j=1,2,3$, are\/
$\Bbb K$-linearly independent and span\/ $\Bbb P_\Bbb KK\subset\Bbb P_\Bbb K(A\otimes_\Bbb KA)$ such that\/
$(L\times\Bbb P_\Bbb KA)\cap\Bbb P_\Bbb KK$ consists only of the points\/ $p_i\otimes p'_i$, $i=1,2,3$.}

\medskip

{\bf Proof.} We can assume that $\sum_ip_i=\sum_jq_j=0$ and $q_2=ap_1+p_2+bq_1$, $a,b\in\Bbb K$. Since neither $p_1,p_2,q_2$ nor
$p_1,p_2,q_3$ lie on a same line, we have $b\ne0,-1$. Suppose that
$$(p_1+cp_2)\otimes d=\sum_ip_i\otimes a_ip'_i+\sum_jq_j\otimes b_jq'_j.$$
Then $p_1\otimes d_1+p_2\otimes d_2+q_1\otimes(b_1q'_1-b_3q'_3+bb_2q'_2-bb_3q'_3)=0$ for suitable $d_1,d_2\in A$. Since $p_1,p_2,q_1$
are not on a same line, we obtain $b_1q'_1+bb_2q'_2-(1+b)b_3q'_3=0$. As $q'_1,q'_2,q'_3$ are not on a same line and $b\ne0,-1$, we
obtain $b_1=b_2=b_3=0$ and
$p_1\otimes d+p_2\otimes cd=\sum_ip_i\otimes a_ip'_i=p_1\otimes(a_1p'_1-a_3p'_3)+p_2\otimes(a_2p'_2-a_3p'_3)$, implying
$d=a_1p'_1-a_3p'_3$ and $cd=a_2p'_2-a_3p'_3$. As $p'_1,p'_2,p'_3$ are not on a same line, we conclude that $a_2=0$. Furthermore,
either $c=a_3=0$ or $c\ne0=a_1$ and $-ca_3=-a_3$. In the first case, we obtain $(p_1+cp_2)\otimes d=p_1\otimes a_1p'_1$. In the second
case, we get $c=1$ and $(p_1+cp_2)\otimes d=p_3\otimes a_3p'_3$
$_\blacksquare$

\medskip

{\bf6.2.~Lemma.} {\sl Let\/ $\varphi_i:D\hookrightarrow\Bbb P_\Bbb KA$, $i=1,2$, be isomorphisms with plane cubics. Then the curve\/
$(\varphi_1\times\varphi_2)D\subset\Bbb P_\Bbb KA\times\Bbb P_\Bbb KA\subset\Bbb P_\Bbb K(A\otimes_\Bbb KA)$ has degree\/ $6$ in\/
$\Bbb P_\Bbb K(A\otimes_\Bbb KA)$ and spans a linear subspace of dimension\/ $\le5$.}

\medskip

{\bf Proof.} Any plane cubic $D$ is a Cohen-Macaulay and Gorenstein scheme with the trivial dualizing sheaf. We obtain the
usual\footnote{I am grateful to Dimitri Markushevich who indicated that the usual Riemann-Roch formula holds for any plane cubic.}
Riemann-Roch formula $h^0\Cal L-h^0\Cal L^\vee=\deg\Cal L$ valid for every line bundle $\Cal L$ over $D$ and deduce that
$h^0\Cal L=\deg\Cal L$ if $\deg\Cal L>0$. As the embedding $\varphi_i:D\hookrightarrow\Bbb P_\Bbb KA$ is given by a line bundle
$\Cal L_i$ of degree $3$, the embedding
$\varphi_1\times\varphi_2:D\hookrightarrow\Bbb P_\Bbb KA\times\Bbb P_\Bbb KA\hookrightarrow P_\Bbb K(A\otimes_\Bbb KA)$ is given by
the line bundle $\Cal L_1\otimes_{\Cal O_D}\Cal L_2$ of degree $6$. Therefore, the image of $\varphi_1\times\varphi_2$ lies in a
linear subspace of $P_\Bbb K(A\otimes_\Bbb KA)$ of dimension $\le5$
$_\blacksquare$

\medskip

{\bf6.3.~Lemma.} {\sl Up to isotopy, there exists a unique\/ $3$-dimensional generic\/ $\Bbb K$-algebra whose zero divisors scheme is
a line plus a double line\/{\rm;} one can take the identity for the isomorphism\/ $\varphi:D_1\to D_2$.}

\medskip

{\bf Proof.} We work in terms of the matrix $\Phi(x_0,x_1,x_2)$ of linear forms in projective coordinates $x_0,x_1,x_2$ on
$\Bbb P_\Bbb KA$ that describe left multiplications by the elements of $A$. After a suitable isotopy, we can assume that $\varphi$
provides the identity on the variety given by $x_0x_1=0$. So, $\Phi(0,x_1,x_2)\left[\smallmatrix0\\x_1\\x_2\endsmallmatrix\right]=0$,
implying that $\Phi(x_0,x_1,x_2)=[*\ f_2x_0+fx_2\ f_3x_0-fx_1]$ for some columns $f_2,f_3,f\in\Bbb K^3$. From
$\Phi(x_0,0,x_2)\left[\smallmatrix x_0\\0\\x_2\endsmallmatrix\right]=0$, we conclude that
$$\Phi(x_0,x_1,x_2)=[f_1x_1-f_3x_2\ f_2x_0+fx_2\ f_3x_0-fx_1]$$
with $f_1,f_2,f_3,f\in\Bbb K^3$. Since
$$x_0^2x_1=\det\Phi(x_0,x_1,x_2)=\det[f_1\ f_2\ f_3]x_0^2x_1-\det[f_1\ f_2\ f]x_0x_1^2+\det[f_1+f_2\ f\ f_3]x_0x_1x_2,$$
we obtain $\det[f_1\ f_2\ f_3]=1$, $\det[f_1\ f_2\ f]=0$, and $\det[f_1+f_2\ f\ f_3]=0$. Acting on $\Phi(x_0,x_1,x_2)$ by
$\GL_3\Bbb K$ from the left, we can assume that $[f_1\ f_2\ f_3]=1$ and infer that $f$ has the form
$f=\left[\smallmatrix k\\k\\0\endsmallmatrix\right]$, $k\in\Bbb K$. Therefore,
$\Phi(x_0,x_1,x_2)=\left[\smallmatrix x_1&kx_2&-kx_1\\0&x_0+kx_2&-kx_1\\-x_2&0&x_0\endsmallmatrix\right]$. Since $\Phi(0,0,1)$ has
rank $2$, we have $k\ne0$. After elementary transformations, we obtain
$\Phi(x_0,x_1,x_2)=\left[\smallmatrix kx_1&kx_2&-kx_1\\-kx_1&x_0&0\\-kx_2&0&x_0\endsmallmatrix\right]$. Denoting by $x_1,x_2$ the
former $kx_1,kx_2$, we get $\Phi(x_0,x_1,x_2)=\left[\smallmatrix x_1&x_2&-x_1\\x_1&x_0&0\\-x_2&0&x_0\endsmallmatrix\right]$. One can
easily see that the isomorphism $\varphi:D_1\to D_2$ is the identity
$_\blacksquare$

\medskip

{\bf6.4.~Lemma.} {\sl Up to isotopy, there exists a unique\/ $3$-dimensional generic\/ $\Bbb K$-algebra whose zero divisors scheme is
a triple line\/{\rm;} one can take the identity for the isomorphism\/ $\varphi:D_1\to D_2$.}

\medskip

{\bf Proof.} After a suitable isotopy, we can assume that $\varphi$ provides the identity on the line of left zero divisors given by
$x_0=0$. As in the proof of Lemma 6.3, we have $\Phi(0,x_1,x_2)\left[\smallmatrix0\\x_1\\x_2\endsmallmatrix\right]=0$, implying that
$\Phi(x_0,x_1,x_2)=[f_1x_0+g_1x_1+g_2x_2\ f_2x_0+g_3x_2\ f_3x_0-g_3x_1]$ for $f_i,g_j\in\Bbb K^3$, $i,j=1,2,3$. Since
$$x_0^3=\det\Phi(x_0,x_1,x_2)=\big(\det[g_1\ f_2\ f_3]-\det[f_1\ f_2\ g_3]\big)x_0^2x_1+\big(\det[g_1\ g_3\ f_3]-\det[g_2\ f_2\
g_3]\big)x_0x_1x_2+$$
$$+\big(\det[f_1\ g_3\ f_3]+\det[g_2\ f_2\ f_3]\big)x_0^2x_2-\det[g_1\ f_2\ g_3]x_0x_1^2+\det[g_2\ g_3\ f_3]x_0x_2^2+\det[f_1\ f_2\
f_3]x_0^3,$$
after acting on $\Phi(x_0,x_1,x_2)$ by $\GL_3\Bbb K$ from the left, we obtain $[f_1\ f_2\ f_3]=1$ and
$$g_{11}=g_{33},\quad g_{11}g_{23}-g_{13}g_{21}=g_{12}g_{33}-g_{13}g_{32},\quad g_{23}+g_{12}=0,\quad g_{11}g_{33}=g_{13}g_{31},\quad
g_{12}g_{23}=g_{13}g_{22},$$
where $[g_{ij}]:=[g_1\ g_2\ g_3]$. This means that
$$[g_1\ g_2\ g_3]=\left[\smallmatrix g_{11}&g_{12}&g_{13}\\g_{21}&g_{22}&-g_{12}\\g_{31}&g_{32}&g_{11}\endsmallmatrix\right],\
g_{11}^2=g_{13}g_{31},\ g_{11}g_{12}+g_{13}g_{21}=g_{13}g_{32}-g_{11}g_{12},\ g_{12}^2+g_{13}g_{22}=0.\leqno{\bold{(6.5)}}$$
If $g_{13}=0$, we have $g_{11}=g_{12}=g_{13}=0$ and
$\Phi(0,x_1,x_2)=\left[\smallmatrix0&0&0\\g_{21}x_1+g_{22}x_2&0&0\\g_{31}x_1+g_{32}x_2&0&0\endsmallmatrix\right]$ has rank $\le1$.
Hence, $g_{13}\ne0$. It follows from (6.5) that
$$[g_1\ 0\ -g_3]\left[\smallmatrix
g_{13}&0&0\\-g_{12}&g_{13}&0\\g_{11}&0&g_{13}\endsmallmatrix\right]=
\left[\smallmatrix0&0&-g_{13}^2\\u^2&0&g_{12}g_{13}\\0&0&-g_{11}g_{13}\endsmallmatrix\right],\qquad[g_2\ g_3\ 0]\left[\smallmatrix
g_{13}&0&0\\-g_{12}&g_{13}&0\\g_{11}&0&g_{13}\endsmallmatrix\right]=
\left[\smallmatrix0&g_{13}^2&0\\0&-g_{12}g_{13}&0\\u^2&g_{11}g_{13}&0\endsmallmatrix\right],$$
where $u^2:=g_{21}g_{13}+g_{12}g_{11}=g_{32}g_{13}-g_{11}g_{12}$. Since
$$\left[\smallmatrix g_{13}&0&0\\g_{12}&g_{13}&0\\-g_{11}&0&g_{13}\endsmallmatrix\right]\left[\smallmatrix0&0&-g_{13}^2\\
u^2&0&g_{12}g_{13}\\0&0&-g_{11}g_{13}\endsmallmatrix\right]=\left[\smallmatrix0&0&-g_{13}^3\\g_{13}u^2&0&0\\0&0&0
\endsmallmatrix\right],\qquad\left[\smallmatrix g_{13}&0&0\\g_{12}&g_{13}&0\\-g_{11}&0&g_{13}\endsmallmatrix\right]
\left[\smallmatrix0&g_{13}^2&0\\0&-g_{12}g_{13}&0\\
u^2&g_{11}g_{13}&0\endsmallmatrix\right]=\left[\smallmatrix0&g_{13}^3&0\\0&0&0\\g_{13}u^2&0&0\endsmallmatrix\right],$$
and
$\Phi(x_0,x_1,x_2)=[f_1\ f_2\ f_3]x_0+[g_1\ 0\ -g_3]x_1+[g_2\ g_3\ 0]x_2$,
we obtain
$$\Phi'(x_0,x_1,x_2):=G^{-1}\Phi(x_0,x_1,x_2)G=\left[\smallmatrix
x_0&g_{13}x_2&-g_{13}x_1\\u^2g_{13}^{-1}x_1&x_0&0\\u^2g_{13}^{-1}x_2&0&x_0\endsmallmatrix\right],$$
where
$G:=\left[\smallmatrix g_{13}&0&0\\-g_{12}&g_{13}&0\\g_{11}&0&g_{13}\endsmallmatrix\right]$. If $u=0$, then $\Phi'(0,x_1,x_2)$ has
rank $\le1$. So, $u\ne0$. Multiplying the first column, second line, and third line by $u^{-1}g_{13}$ in $\Phi'(x_0,x_1,x_2)$, we
arrive at
$\left[\smallmatrix
u^{-1}g_{13}x_0&g_{13}x_2&-g_{13}x_1\\g_{13}x_1&u^{-1}g_{13}x_0&0\\g_{13}x_2&0&u^{-1}g_{13}x_0\endsmallmatrix\right]$.
Denoting the former $u^{-1}g_{13}x_0,g_{13}x_1,g_{13}x_2$ by $x_0,x_1,x_2$, we obtain
$\left[\smallmatrix x_0&x_2&-x_1\\x_1&x_0&0\\x_2&0&x_0\endsmallmatrix\right]$. One can easily see that the isomorphism
$\varphi:D_1\to D_2$ is the identity
$_\blacksquare$

\medskip

{\bf Proof of Theorem 1.3.2.} Up to isotopy, any algebra $A$ can be described by
$\Bbb P_\Bbb KK\subset\Bbb P_\Bbb K(A\otimes_\Bbb KA)$, where $K$ stands for the kernel of the multiplication
$A\otimes_\Bbb KA\overset\cdot\to\to A$. In view of the Segre embedding
$\Bbb P_\Bbb KA\times\Bbb P_\Bbb KA\subset\Bbb P_\Bbb K(A\otimes_\Bbb KA)$, the scheme $D$ of the pairs of zero divisors equals
$D=(\Bbb P_\Bbb KA\times\Bbb P_\Bbb KA)\cap\Bbb P_\Bbb KK$.

One can easily see that the Lie algebra $A:=\ssl_2\Bbb K$ is a generic algebra with $D_1=\Bbb P_\Bbb KA$; that the scheme $D$
coincides with the diagonal $\Delta\subset\Bbb P_\Bbb KA\times\Bbb P_\Bbb KA$; and that
$\Delta=(\Bbb P_\Bbb KA\times\Bbb P_\Bbb KA)\cap\Bbb P_\Bbb KK$ spans $\Bbb P_\Bbb KK$, where $K=\Sym^2A\subset A\otimes_\Bbb KA$ is
the symmetric square of $A$, $\dim_\Bbb KK=6$. Hence, up to isotopy, $\ssl_2\Bbb K$ is a unique generic algebra with
$D_1=\Bbb P_\Bbb KA$.

It follows that $\dim_\Bbb KK=6$ for any generic algebra. Indeed, if $\dim_\Bbb KK>6$, then $D_1=\Bbb P_\Bbb KA$ and $A$ has to be
isotopic to $\ssl_2\Bbb K$, a contradiction.

Let $A$ be a generic $3$-dimensional algebra whose $D_1$ is a reduced cubic, i.e., a cubic without multiple components.

If the isomorphism $\varphi:D_1\to D_2$ is projective, we can assume that $D=\Delta_{D_1}:=(D_1\times D_1)\cap\Delta$. A~simple
straightforward verification shows that $\Delta_{D_1}$ spans $\Bbb P_\Bbb K\Sym^2A$ (for a smooth $D_1$, one may use the Riemann-Roch
theorem as above). In other words, we arrive again at $\ssl_2\Bbb K$. A contradiction.

So, we assume that $\varphi$ is not projective. This implies that, for any $3$ distinct points
$p_1,p_2,p_3\in D_1\subset\Bbb P_\Bbb KA$ lying on a same line transversal to $D_1$, the points
$\varphi p_1,\varphi p_2,\varphi p_3\in D_2\subset\Bbb P_\Bbb KA$ do not lie on a same line. We can choose points $p_i,q_j\in D_1$
with $p'_i:=\varphi p_i$ and $q'_j:=\varphi q_i$ subject to the conditions of Lemma~6.1 and with $L,M$ transversal to $D_1$. By
Lemmas 6.1 and 6.2, $\Bbb P_\Bbb KK\subset\Bbb P_\Bbb K(A\otimes_\Bbb KA)$ is spanned by the pairs $(d_1,\varphi d_1)$, $d_1\in D_1$.

Conversely, let us be given a nonprojective isomorphism $\varphi:D_1\to D_2$ between reduced cubics. By Lemmas 6.1 and 6.2, the curve
$D':=\big\{(d_1,\varphi d_1)\mid d_1\in D_1\big\}$ of degree $6$ in $\Bbb P_\Bbb K(A\otimes_\Bbb KA)$ spans a $5$-dimensional linear
subspace $\Bbb P_\Bbb KK\subset\Bbb P_\Bbb K(A\otimes_\Bbb KA)$.

For any smooth point $p_1\in D_1$, there is a line $L\ni p_1$ transversal to $D_1$. Taking a generic line $M$ transversal to $D_1$, we
conclude from Lemma 6.1 that $(L\times\Bbb P_\Bbb KA)\cap\Bbb P_\Bbb KK$ consists of $3$ points. In particular,
$(p_1\times\Bbb P_\Bbb KA)\cap\Bbb P_\Bbb KK=(p_1,\varphi p_1)$ and there is a point in $L$ that is not a left zero divisor of the
algebra $A$ given by $K$. This implies that $D_1$ is the scheme of left zero divisors of $A$. By symmetry, $D_2$ is the scheme of
right zero divisors of $A$.

We need to show that $D'=D:=(\Bbb P_\Bbb KA\times\Bbb P_\Bbb KA)\cap\Bbb P_\Bbb KK$. Since $\Bbb P_\Bbb KA\times\Bbb P_\Bbb KA$ has
degree $6$ in $\Bbb P_\Bbb K(A\otimes_\Bbb KA)$, we can assume that $\dim D\ge2$. So, $D\supset C_1\times C_2$, where $C_i$ is a
component of $D_i$. As there is a point $p_1\in C_1$ that is smooth in $D_1$ and $p_1\times C_2\subset D\subset\Bbb P_\Bbb KK$, we
arrive at a contradiction.

The remaining case of $D_1$ with multiple components is considered in Lemmas 6.3 and 6.4
$_\blacksquare$

\medskip

When studying noncommutative projective planes, A.~Bondal and A.~Polishchuk [BoP] classified the so-called geometric tensors. This
classification is almost equivalent to that of generic algebras (see [BoP, Table, p.~36] for details). The algebras given by matrices
$\left[\smallmatrix x_0&0&0\\-ux_2&x_0&ux_1\\-u^{-1}x_1&u^{-1}x_2&x_0\endsmallmatrix\right]$ of linear forms, where $0\ne u\in\Bbb K$,
constitute the difference between generic algebras and geometric tensors. In this case, each of $D_1,D_2$ is a conic plus a line and
the isomorphism $\varphi:D_1\to D_2$ maps the line to the conic and the conic to the line.

\bigskip

\centerline{\bf References}

\medskip

[AGG] S.~Anan$'$in, C.~H.~Grossi, N.~Gusevskii, {\it Complex hyperbolic structures on disc bundles over surfaces,}
Int.~Math.~Res.~Not.~{\bf2011} (2011), no.~19, 4295--4375, see also http://arxiv.org/abs/math/0511741

[AGr] S.~Anan$'$in, C.~H.~Grossi, {\it Coordinate-free classic geometries,} Mosc.~Math.~J.~{\bf11} (2011), 633--655, see also
http://arxiv.org/abs/math/0702714

[AGS] S.~Anan$'$in, C.~H.~Grossi, J.~S.~S.~da Silva, {\it Poincare's polyhedron theorem for cocompact groups in dimension\/ $4$,}
accepted in Mosc.~Math.~J.~(2014), see also http://arxiv.org/abs/1112.5740

[Ana] S.~Anan$'$in, {\it Research plan for\/} 2013--2014

[BoP] A.~Bondal, A.~E.~Polishchuk, {\it Homological properties of associative algebras\/{\rm:} the method of helices,}
Russian Acad.~Sci.~Izv.~Math.~42:2 (1994), 219--260

[Gol] W.~M.~Goldman, {\it Complex hyperbolic geometry,} Oxford Mathematical Monographs. Oxford Science Publications. The Clarendon
Press, Oxford University Press, New York, 1999, xx+316 pp.

\enddocument